\documentclass[11pt,a4paper]{amsart} 
\usepackage[all]{xy}
\SelectTips{cm}{}
\calclayout \makeatletter
\def\serieslogo@{} 
\def\@setcopyright{} 
\makeatother
\title[Cohomological quotients]{
Cohomological quotients and smashing localizations
} 
\author[Henning Krause]{Henning Krause}
\address{Henning Krause\\ Institut f\"ur Mathematik\\
Universit\"at Paderborn\\ 33095 Paderborn\\ Germany.}
\email{hkrause@math.upb.de}
\thanks{Version from \today.}  
\newtheorem{lem}{Lemma}[section]
\newtheorem{prop}[lem]{Proposition}
\newtheorem{cor}[lem]{Corollary}
\newtheorem{thm}[lem]{Theorem}

\newtheorem*{Thm1}{Theorem~1}
\newtheorem*{Thm2}{Theorem~2}
\newtheorem*{Thm3}{Theorem~3}

\newtheorem*{Cor}{Corollary}
\newtheorem*{Prop}{Proposition}
\newtheorem*{telconj}{Telescope Conjecture}

\theoremstyle{remark}
\newtheorem{rem}[lem]{Remark}

\theoremstyle{definition}
\newtheorem{exm}[lem]{Example}
\newtheorem{defn}[lem]{Definition}
\newtheorem*{Defn}{Definition}
\numberwithin{equation}{section}

\renewcommand{\mod}{\operatorname{mod}\nolimits}
\newcommand{\colim}{\operatorname{colim}\nolimits}
\newcommand{\Ann}{\operatorname{Ann}\nolimits}
\newcommand{\proj}{\operatorname{proj}\nolimits}
\newcommand{\rad}{\operatorname{rad}\nolimits}

\newcommand{\id}{\operatorname{id}\nolimits}

\newcommand{\Mod}{\operatorname{Mod}\nolimits}

\newcommand{\Hom}{\operatorname{Hom}\nolimits}
\newcommand{\Ho}{\operatorname{Ho}\nolimits}
\newcommand{\HOM}{\operatorname{\mathcal H \mathit o \mathit m}\nolimits}
\renewcommand{\Im}{\operatorname{Im}\nolimits}
\newcommand{\Ker}{\operatorname{Ker}\nolimits}
\newcommand{\Coker}{\operatorname{Coker}\nolimits}
\newcommand{\Cone}{\operatorname{Cone}\nolimits}

\newcommand{\Ext}{\operatorname{Ext}\nolimits}
\newcommand{\Tor}{\operatorname{Tor}\nolimits}
\newcommand{\Filt}{\operatorname{Filt}\nolimits}

\newcommand{\umod}{\operatorname{\underline{mod}}\nolimits}

\newcommand{\uHom}{\operatorname{\underline{Hom}}\nolimits}

\newcommand{\Ab}{\mathrm{Ab}}
\newcommand{\op}{\mathrm{op}}
\newcommand{\Id}{\mathrm{Id}}
\newcommand{\ex}{\mathrm{ex}}

\newcommand{\coh}{\mathrm{coh}}
\newcommand{\comp}{\mathop{\raisebox{+.3ex}{\hbox{$\scriptstyle\circ$}}}}
\newcommand{\lto}{\longrightarrow}

\def\li{\varinjlim}

\def\a{\alpha}
\def\b{\beta}
\def\e{\varepsilon}
\def\d{\delta}
\def\g{\gamma}
\def\p{\phi}

\def\k{\kappa}

\def\Ga{\Gamma}
\def\La{\Lambda}
\def\Si{\Sigma}

\def\A{{\mathcal A}}
\def\B{{\mathcal B}}
\def\C{{\mathcal C}}
\def\D{{\mathcal D}}
\def\E{{\mathcal E}}
\def\I{{\mathcal I}}

\def\F{{\mathcal F}}

\def\K{{\mathcal K}}
\def\L{{\mathcal L}}
\def\M{{\mathcal M}}

\def\R{{\mathcal R}}
\def\S{{\mathcal S}}
\def\X{{\mathcal X}}

\def\T{{\mathcal T}}

\begin{document}
\begin{abstract}
The quotient of a triangulated category modulo a subcategory was
defined by Verdier. Motivated by the failure of the telescope
conjecture, we introduce a new type of quotients for any triangulated
category which generalizes Verdier's construction.  Slightly
simplifying this concept, the cohomological quotients are flat
epimorphisms, whereas the Verdier quotients are Ore localizations. For
any compactly generated triangulated category $\S$, a bijective
correspondence between the smashing localizations of $\S$ and the
cohomological quotients of the category of compact objects in $\S$ is
established. We discuss some applications of this theory, for instance
the problem of lifting chain complexes along a ring homomorphism. This
is motivated by some consequences in algebraic $K$-theory and
demonstrates the relevance of the telescope conjecture for derived
categories. Another application leads to a derived analogue of an
almost module category in the sense of Gabber-Ramero.  It is shown
that the derived category of an almost ring is of this form.
\end{abstract}
\maketitle
\tableofcontents
\section*{Introduction}
The telescope conjecure from stable homotopy theory is a fascinating
challenge for topologists and algebraists.  It is a conjecture about
smashing localizations, saying roughly that every smashing
localization is a finite localization.  The failure of this conjecture
forces us to develop a general theory of smashing localizations which
covers the ones which are not finite. This is precisely the subject of
the first part of this paper. The second part discusses some
applications of the general theory in the context of derived
categories of associative rings. In fact, we demonstrate the relevance
of the telescope conjecture for derived categories, by studying some
applications in algebraic $K$-theory and in almost ring theory.

Let us describe the main concepts and results from this paper.  We fix
a compactly generated triangulated category $\S$, for example, the
stable homotopy category of CW-spectra or the unbounded derived
category of an associative ring. A {\em smashing localization functor}
is by definition an exact functor $F\colon \S\to\T$ between
triangulated categories having a right adjoint $G$ which preserves all
coproducts and satisfies $F\comp G\cong\Id_\T$. Such a functor induces
an exact functor $F_c\colon \S_c\to\T_c$ between the full
subcategories of compact objects, and the telescope conjecture
\cite{B,R} claims that the induced functor $\S_c/{\Ker F_c}\to\T_c$ is
an equivalence up to direct factors. Here, $\Ker F_c$ denotes the full
triangulated subcategory of objects $X$ in $\S_c$ such that $F_cX=0$,
and $\S_c/{\Ker F_c}$ is the quotient in the sense of Verdier
\cite{V}. The failure of the telescope conjecture \cite{Ke,MRS}
motivates the following generalization of Verdier's definition of a
quotient of a triangulated category. To be precise, there are examples
of proper smashing localization functors $F$ where $\Ker
F_c=0$. Nonetheless, the functor $F_c$ is a cohomological quotient
functor in the following sense.

\begin{Defn}
Let $F\colon\C\to\D$ be an exact functor between triangulated
categories. We call $F$ a {\em cohomological quotient functor} if for
every cohomological functor $H\colon\C\to \A$ satisfying $\Ann
F\subseteq\Ann H$, there exists, up to a unique isomorphism, a unique
cohomological functor $H'\colon\D\to \A$ such that $H=H'\comp F$.
\end{Defn}
Here, $\Ann F$ denotes the ideal of all maps $\p$ in $\C$ such that
$F\p=0$. The property of $F$ to be a cohomological quotient functor
can be expressed in many ways, for instance more elementary as
follows: every object in $\D$ is a direct factor of some object in the
image of $F$, and every map $\a\colon FX\to FY$ in $\D$ can be
composed with a split epimorphism $F\pi\colon FX'\to FX$ such that
$\a\comp F\pi$ belongs to the image of $F$.

Our main result shows a close relation between cohomological quotient
functors and smashing localizations.

\begin{Thm1}
Let $\S$ be a compactly generated triangulated category, and let 
$F\colon\S_c\to\T$ be a cohomological quotient functor. 
Denote by $\R$ the full subcategory of objects $X$ in $\S$ such that every
map $C\to X$ from a compact object $C$ factors through some map in $\Ann F$.
\begin{enumerate}
\item The category $\R$ is a triangulated subcategory of $\S$ and the
quotient functor $\S\to\S/\R$ is a smashing localization functor which
induces a fully faithful and exact functor $\T\to\S/\R$ making the
following diagram commutative.
$$\xymatrix{\S_c\ar[rr]^F\ar[d]^{\mathrm{inc}}&&\T\ar[d]\\
\S\ar[rr]^{\mathrm{can}}&&\S/\R}$$
\item The triangulated category $\S/\R$ is compactly generated and the
subcategory of compact objects is precisely the closure of the image
of $\T\to\S/\R$ under forming direct factors.
\item There exists a fully faithful and exact functor
$G\colon\T\to\S$ such that
$$\S(X,GY)\cong\T(FX,Y)$$
for all $X$ in $\S_c$ and $Y$ in $\T$.
\end{enumerate}
\end{Thm1}

One may think of this result as a generalization of the
localization theorem of Neeman-Ravenel-Thomason-Trobaugh-Yao
\cite{N2,R,TT,Y}.  To be precise, Neeman et al.\ considered
cohomological quotient functors of the form $\S_c\to\S_c/\R_0$ for
some triangulated subcategory $\R_0$ of $\S_c$ and analyzed the
smashing localization functor $\S\to\S/\R$ where $\R$ denotes the
localizing subcategory generated by $\R_0$.

Our theorem provides a bijective correspondence between smashing
localizations of $\S$ and cohomological quotients of $\S_c$; it
improves a similar correspondence \cite{K} -- the new ingredient in
our proof being a recent variant \cite{K2} of Brown's
Representability Theorem \cite{Br}. The essential invariant of a
cohomological quotient functor $F\colon\S_c\to\T$ is the ideal $\Ann
F$. The ideals of $\S_c$ which are of this form are called {\em exact}
and are precisely those satisfying the following properties:
\begin{enumerate}
\item $\mathfrak I^2=\mathfrak I$.
\item $\mathfrak I$ is {\em saturated}, that is, for every exact triangle
$X'\stackrel{\a}\to X\stackrel{\b}\to X''\to\Si X'$ and every map
$\p\colon X\to Y$ in $\S_c$, we have that $\p\comp\a,\b\in\mathfrak
I$ implies $\p\in\mathfrak I$.
\item $\Si\mathfrak I=\mathfrak I$.
\end{enumerate}

Let us rephrase the telescope conjecture in terms of exact ideals and
cohomological quotient functors. To this end, recall that a
subcategory of $\S$ is {\em smashing} if it is of the form $\Ker F$
for some smashing localization functor $F\colon\S\to\T$.

\begin{Cor} 
The telescope conjecture for $\S$ is equivalent to each of the
following statements.
\begin{enumerate}
\item Every smashing subcategory of $\S$ is generated by compact objects.
\item Every exact ideal is generated by idempotent elements.
\item Every cohomological quotient functor $F\colon\S_c\to\T$ induces
up to direct factors an equivalence $\S_c/{\Ker F}\to\T$.
\item Every flat epimorphism $F\colon \S_c\to\T$ satisfying $\Si(\Ann
F)=\Ann F$ is an Ore localization.
\end{enumerate}
\end{Cor}
This reformulation of the telescope conjecture is based on our
approach to view a triangulated category as a ring with several
object. In this setting, the cohomological quotient functors are the
flat epimorphisms, whereas the Verdier quotient functors are the Ore
localizations. The reformulation in terms of exact ideals refers to
the classical problem from ring theory of finding idempotent
generators for an idempotent ideal, studied for instance by Kaplansky
\cite{J} and Auslander \cite{A}. We note that the telescope conjecture
becomes a statement about the category of compact objects. Moreover,
we see that the smashing subcategories of $\S$ form a complete lattice
which is isomorphic to the lattice of exact ideals in $\S_c$.

The second part of this paper is devoted to studying non-commutative
localizations of rings. We do this by using unbounded derived
categories and demonstrate that the telescope conjecture is relevant
in this context. This is inspired by recent work of Neeman and Ranicki
\cite{NR}.  They study the problem of lifting chain complexes up to
homotopy along a ring homomorphism $R\to S$. To make this
precise, let us denote by $\mathbf K^b(R)$ the homotopy category of
bounded complexes of finitely generated projective $R$-modules.
\begin{enumerate}
\item We say that the {\em chain complex lifting problem} has a
positive solution, if every complex $Y$ in $\mathbf K^b(S)$ such that
for each $i$ we have $Y^i=P^i\otimes_RS$ for some finitely generated
projective $R$-module $P^i$, is isomorphic to $X\otimes_RS$ for some
complex $X$ in $\mathbf K^b(R)$.
\item We say that the {\em chain map lifting problem} has a positive
solution, if for every pair $X,Y$ of complexes in $\mathbf K^b(R)$ and
every map $\a\colon X\otimes_RS\to Y\otimes_RS$ in $\mathbf K^b(S)$,
there are maps $\phi\colon X'\to X$ and $\a'\colon X'\to Y$ in
$\mathbf K^b(R)$ such that $\phi\otimes_RS$ is invertible and
$\a=\a'\otimes_RS\comp(\phi\otimes_RS)^{-1}$ in $\mathbf K^b(S)$.
\end{enumerate}
Note that complexes can be lifted whenever maps can be lifted. For
example, maps and complexes can be lifted if $R\to S$ is a commutative
localization. However, there are obstructions in the non-commutative
case, and this leads to the concept of a homological
epimorphism. Recall from \cite{GL} that $R\to S$ is a {\em homological
epimorphism} if $S\otimes_RS\cong S$ and $\Tor_i^R(S,S)=0$ for all
$i\geq 1$. For example, every commutative localization is a flat
epimorphism and therefore a homological epimorphism.  The following
observation is crucial for both lifting problems.

\begin{Prop} 
A ring homomorphism $R\to S$ is a homological epimorphism if and only
if $-\otimes_RS\colon\mathbf K^b(R)\to\mathbf K^b(S)$ is a
cohomological quotient functor.
\end{Prop}

This shows that we can apply our theory of cohomological quotient
functors, and we see that the telescope conjecture for the unbounded
derived category $\mathbf D(R)$ of a ring $R$ becomes relevant. In
particular, we obtain a non-commutative analogue of Thomason-Trobaugh's
localization theorem for algebraic $K$-theory \cite{TT}.

\begin{Thm2}
Let $R$ be a ring such that the telescope conjecture holds true for
$\mathbf D(R)$. Then the chain map lifting problem has a positive
solution for a ring homomorphism $f\colon R\to S$ if and only if $f$
is a homological epimorphism. Moreover, in this case $f$ induces a
sequence
$$K(R,f)\lto K(R)\lto K(S)$$ of $K$-theory spectra which is a homotopy
fibre sequence, up to failure of surjectivity of $K_0(R)\to K_0(S)$.
In particular, there is induced a long exact sequence
$$\cdots \lto K_1(R)\lto K_1(S)\lto K_0(R,f)\lto K_0(R)\lto K_0(S)$$
of algebraic $K$-groups.
\end{Thm2}

Unfortunately, not much seems to be known about the telescope
conjecture for derived categories.  Note that the telescope conjecture
has been verified for $\mathbf D(R)$ provided $R$ is commutative
noetherian \cite{N}. On the other hand, there are counter examples
which arise from homological epimorphisms where not all chain maps can
be be lifted \cite{Ke}.

In the final part of this paper, we introduce the derived analogue of
an almost module category in the sense of \cite{GR}. In fact, there is
a striking parallel between almost rings and smashing localizations:
both concepts depend on an idempotent ideal. Given a ring $R$ and an
idempotent ideal $\mathfrak a$, the category of {\em almost modules} is
by definition the quotient
$$\Mod (R,\mathfrak a)=\Mod R/(\mathfrak a^\perp),$$ where $\Mod R$
denotes the category of right $R$-modules and $\mathfrak a^\perp$
denotes the Serre subcategory of $R$-modules annihilated by $\mathfrak
a$. Given an idempotent ideal $\mathfrak I$ of $\mathbf K^b(R)$ which
satifies $\Si\mathfrak I=\mathfrak I$, the objects in $\mathbf D(R)$
which are annihilated by $\mathfrak I$ form a triangulated
subcategory, and we call the quotient category
$$\mathbf D(R,\mathfrak I)=\mathbf D(R)/({\mathfrak I}^\perp)$$ an
{\em almost derived category}. It turns out that the almost derived
categories are, up to equivalence, precisely the smashing
subcategories of $\mathbf D(R)$. Moreover, as one should expect, the
derived category of an almost ring is an almost derived category.

\begin{Thm3}
Let $R$ be a ring and $\mathfrak a$ be an idempotent ideal such that
$\mathfrak a\otimes_R\mathfrak a$ is flat as left $R$-module.  Then
the maps in $\mathbf K^b(R)$ which annihilate all suspensions of the
mapping cone of the natural map $\mathfrak a\otimes_R\mathfrak a\to R$
form an idempotent ideal $\mathfrak A$, and $\mathbf D(R,\mathfrak A)$
is equivalent to the unbounded derived category of $\Mod (R,\mathfrak
a)$.
\end{Thm3}

\subsection*{Acknowledgements} 
I would like to thank Ragnar Buchweitz, Bernhard Keller, and Amnon
Neeman for several stimulating discussions during a visit to the
Mathematical Sciences Institute in Canberra in July 2003.  In
addition, I am grateful to an anonymous referee for a number of
helpful comments.

\section{Modules}

The homological properties of an additive category $\C$ are reflected
by properties of functors from $\C$ to various abelian categories.  In
this context, the abelian category $\Ab$ of abelian groups plays a
special role, and this leads to the concept of a $\C$-module.  In this
section we give definitions and fix some terminology.

Let $\C$ and $\D$ be additive categories. We denote by $\HOM(\C,\D)$
the category of functors from $\C$ to $\D$. The natural
transformations between two functors form the morphisms in this
category, but in general they do not form a set. A category will be
called {\em large} to point out that the morphisms between fixed
objects are not assumed to form a set.

A {\em $\C$-module} is by definition an additive functor
$\C^\op\to\Ab$ into the category $\Ab$ of abelian groups, and we
denote for $\C$-modules $M$ and $N$ by $\Hom_\C(M,N)$ the class of
natural transformations $M\to N$. We write $\Mod\C$ for the category
of $\C$-modules which is large, unless $\C$ is {\em small}, that is,
the isomorphism classes of objects in $\C$ form a set. Note that
$\Mod\C$ is an abelian category. A sequence $L\to M\to N$ of maps
between $\C$-modules is {\em exact} if the sequence $LX\to MX\to NX$
is exact for all $X$ in $\C$.  We denote for every $X$ in $\C$ by
$H_X=\C(-,X)$ the corresponding representable functor and recall that
$\Hom_\C(H_X,M)\cong MX$ for every module $M$ by Yoneda's lemma. It
follows that $H_X$ is a projective object in $\Mod\C$.

A $\C$-module $M$ is called {\em finitely presented} if it fits into
an exact sequence
$$\C(-,X)\lto\C(-,Y)\lto M\to 0$$ with $X$ and $Y$ in $\C$. Note that
$\Hom_\C(M,N)$ is a set for every finitely presented $\C$-module $M$
by Yoneda's lemma. The finitely presented $\C$-modules form an
additive category with cokernels which we denote by $\mod\C$.

Now let $F\colon\C\to\D$ be an additive functor.  This induces the
{\em restriction functor}
$$F_*\colon\Mod\D\lto\Mod\C,\quad M \longmapsto M\comp F,$$
and its left adjoint
$$F^*\colon\Mod\C\lto\Mod\D$$ 
which  sends a $\C$-module $M$, written as a colimit
$M=\colim_{\a\in MX}\C(-,X)$ of representable functors, to
$$F^*M=\colim_{\a\in MX}\D(-,FX).$$ Note that every $\C$-module can be
written as a small colimit of representable functors provided $\C$ is
small.  The finitely presented $\C$-modules are precisely the finite
colimits of representable functors. We denote the restriction of $F^*$ by
$$F^\star\colon\mod\C\lto\mod\D$$ and observe that $F^\star$ is the unique
right exact functor $\mod\C\to\mod\D$ sending $\C(-,X)$ to $\D(-,FX)$
for all $X$ in $\C$.

Finally, we define
\begin{align*}
\Ann F&= \mbox{the ideal of all maps $\p\in\C$ with $F\p=0$, and } \\
\Ker F&= \mbox{the full subcategory of all objects $X\in\C$ with $FX=0$.}
\end{align*}
Recall that an {\em ideal} $\mathfrak I$ in $\C$ consists of subgroups
$\mathfrak I(X,Y)$ in $\C(X,Y)$ for every pair of objects $X,Y$ in
$\C$ such that for all $\p$ in $\mathfrak I(X,Y)$ and all maps
$\a\colon X'\to X$ and $\b\colon Y\to Y'$ in $\C$ the composition
$\b\comp \p\comp \a$ belongs to $\mathfrak I(X',Y')$.  Note that all
ideals in $\C$ are of the form $\Ann F$ for some additive functor
$F$.

Given any class $\Phi$ of maps in $\C$, we say that an object $X$ in
$\C$ is {\em annihilated} by $\Phi$, if $\Phi\subseteq\Ann\C(-,X)$. We
denote by $\Phi^\perp$ the full subcategory of objects in $\C$ which
are annihilated by $\Phi$.

\section{Cohomological functors and ideals}

Let $\C$ be an additive category and suppose $\mod\C$ is abelian.
Note that $\mod\C$ is abelian if and only if every map $Y\to Z$ in
$\C$ has a {\em weak kernel} $X\to Y$, that is, the sequence
$\C(-,X)\to\C(-,Y)\to\C(-,Z)$ is exact. In particular, $\mod\C$ is
abelian if $\C$ is triangulated.  A functor $F\colon\C\to\A$ to an
abelian category $\A$ is called {\em cohomological} if it sends every
weak kernel sequence $X\to Y\to Z$ in $\C$ to an exact sequence $FX\to
FY\to FZ$ in $\A$.  If $\C$ is a triangulated category, then a functor
$F\colon\C\to\A$ is cohomological if and only if $F$ sends every exact
triangle $X\to Y\to Z\to\Si X$ in $\C$ to an exact sequence $FX\to
FY\to FZ\to F\Si X$ in $\A$. The Yoneda functor
$$H_\C\colon\C\longrightarrow\mod\C, \quad X\mapsto H_X=\C(-,X)$$ is
the universal cohomological functor for $\C$. More precisely, for
every abelian category $\A$, the functor
$$\HOM(H_\C,\A)\colon\HOM(\mod\C,\A)\lto\HOM(\C,\A)$$
induces an equivalence
$$\HOM_\ex(\mod\C,\A)\lto\HOM_\coh(\C,\A),$$ where the subscripts
$\ex$ = exact and $\coh$ = cohomological refer to the appropriate full
subcategories; see \cite{F,V} and also \cite[Lemma~2.1]{K}.

Following \cite{K}, we call an ideal $\mathfrak I$ in $\C$ {\em
cohomological} if there exists a cohomological functor $F\colon\C\to
\A$ such that $\mathfrak I=\Ann F$.  For example, if $F\colon\C\to\D$
is an exact functor between triangulated categories, then $\Ann F$ is
cohomological because $\Ann F=\Ann (H_\D\comp F)$.  Note that the
cohomological ideals of $\C$ form a complete lattice,
provided $\C$ is small. For instance, given a family $(\mathfrak
I_i)_{i\in\La}$ of cohomological ideals, we have
$$\inf\mathfrak I_i=\bigcap_i\mathfrak I_i,$$ 
because  $\bigcap_i\mathfrak
I_i=\Ann F$ for
$$F\colon\C\lto\prod_i\A_i, \quad X\mapsto (F_iX)_{i\in\La}$$ where
each $F_i\colon\C\to\A_i$ is a cohomological functor satisfying
$\mathfrak I_i=\Ann F_i$. We obtain $\sup\mathfrak I_i$ by taking the
infimum of all cohomological ideals $\mathfrak J$ with $\mathfrak
I_i\subseteq \mathfrak J$ for all $i\in\La$.

\section{Flat epimorphisms}

The concept of a flat epimorphisms generalizes the classical notion of
an Ore localization.  We study flat epimorphisms of additive
categories, following the idea that an additive category may be viewed
as a ring with several objects. Given a flat epimorphism $\C\to\D$, it
is shown that the maps in $\D$ are obtained from those in $\C$ by a
generalized calculus of fractions. There is a close link between flat
epimorphisms and quotients of abelian categories.  It is the aim of
this section to explain this connection which is summarized in
Theorem~\ref{th:flatepi}.  We start with a brief discussion of
quotients of abelian categories.

Let $\C$ be an abelian category. A full subcategory $\B$ of $\C$ is
called a {\em Serre subcategory} provided that for every exact
sequence $0\to X'\to X\to X'' \to 0$ in $\C$, the object $X$ belongs
to $\B$ if and only if $X'$ and $X''$ belong to $\B$.  The {\em
quotient} $\C/\B$ with respect to a Serre subcategory $\B$ is by
definition the localization $\C[\Phi^{-1}]$, where $\Phi$ denotes the
class of maps $\p$ in $\C$ such that $\Ker\p$ and $\Coker\p$ belong to
$\B$; see \cite{G,GZ}.  The localization functor $Q\colon\C\to\C/\B$
yields for every category $\E$ a functor
$$\HOM(Q,\E)\colon\HOM(\C/\B,\E)\lto\HOM(\C,\E)$$ which induces an
isomorphism onto the full subcategory of functors $F\colon\C\to\E$
such that $F\p$ is invertible for all $\p\in\Phi$. Note that $\C/\B$
is abelian and $Q$ is exact with $\Ker Q=\B$. Up to an equivalence, a
localization functor can be characterized as follows.

\begin{lem}\label{le:abelian1}
Let $F\colon\C\to\D$ be an exact functor between abelian categories.
Then the following are equivalent.
\begin{enumerate}
\item $F$ induces an equivalence $\C/{\Ker F}\to\D$.
\item For every abelian category $\A$, the functor
$$\HOM(F,\A)\colon\HOM_\ex(\D,\A)\lto\HOM_\ex(\C,\A)$$ induces an
equivalence onto the full subcategory of exact functors $G\colon\C\to\A$ satisfying
$\Ker F\subseteq\Ker G$.
\end{enumerate}
\end{lem}
\begin{proof} See \cite[III.1]{G}.
\end{proof}

An exact functor between abelian categories satisfying the equivalent
conditions of Lemma~\ref{le:abelian1} is called an {\em exact quotient
functor}. There is a further characterization in case the functor has
a right adjoint.

\begin{lem}\label{le:abelian2}
Let $F\colon\C\to\D$ be an exact functor between abelian categories
and suppose there is a right adjoint $G\colon\D\to\C$. Then $F$ is a
quotient functor if and only if $G$ is fully faithful. In this case, $G$
identifies $\D$ with the full subcategory of objects $X$ in $\C$
satisfying $\C(\Ker F,X)=0$ and $\Ext_\C^1(\Ker F,X)=0$.
\end{lem}
\begin{proof} See Proposition~III.3 and Proposition~III.5 in \cite{G}.
\end{proof}

Next we analyze an additive functor $F\colon\C\to\D$ in terms of the
induced functor $F^\star\colon\mod\C\to\mod\D$.

\begin{lem}\label{le:maps}
  Let $F\colon\C\to\D$ be an additive functor between additive
  categories and suppose $F^\star\colon\mod\C\to\mod\D$ is an exact
  quotient functor of abelian categories.
\begin{enumerate}
\item Every object in $\D$ is a direct factor of some object in the image of $F$.
\item For every map $\a\colon FX\to FY$ in $\D$, there are maps
$\a'\colon X'\to Y$ and $\pi\colon X'\to X$ in $\C$ such that
$F\a'=\a\comp F\pi$ and $F\pi$ is a split epimorphism.
\end{enumerate}
\end{lem}
\begin{proof}
The functor $F^\star\colon\mod\C\to\mod\D$ is, up to an equivalence, a
localization functor.  Therefore the objects in $\mod\D$ coincide, up
to isomorphism, with the objects in $\mod\C$.  Moreover, the maps in
$\mod\D$ are obtained via a calculus of fractions from the maps in
$\mod\C$; see \cite[I.2.5]{GZ}.

(1) Fix an object $Y$ in $\D$. Then $\D(-,Y)\cong F^\star M$ for some $M$ in
$\mod\C$. If $M$ is a quotient of $\C(-,X)$, then $F^\star M$ is a quotient
of $\D(-,FX)$. Thus $Y$ is a direct factor of $FX$.

(2) Fix a map $\a\colon FX\to FY$.  The corresponding map $\D(-,\a)$
in $\mod\D$ is a fraction, that is, of the form
$$\D(-,FX)=F^\star \C(-,X)\xrightarrow{(F^\star \sigma)^{-1}} F^\star
M\xrightarrow{F^\star \p} F^\star \C(-,Y)=\D(-,FY)$$
for some $M$ in
$\mod\C$; see \cite[I.2.5]{GZ}. Choose an epimorphism $\rho\colon
\C(-,X')\to M$ for some $X'$ in $\C$. Now define $\a'\colon X'\to Y$
by $\C(-,\a')=\p\comp\rho$, and define $\pi \colon X'\to X$ by
$\C(-,\pi)=\sigma\comp\rho$.  Clearly, $F\pi$ is a split epimorphism
since $\D(-,FX)$ is a projective object in $\mod\D$.
\end{proof}

\begin{rem}\label{re:maps}
Conditions (1) and (2) in Lemma~\ref{le:maps} imply that every map
in $\D$ is a direct factor of some map in the image of $F$. To be
precise, we say that a map $\a\colon X\to X'$ is a {\em direct factor}
of a map $\b\colon Y\to Y'$ if there is a commutative diagram
$$\xymatrix{ X\ar[d]^{\a}\ar[r]^{\e}&Y\ar[d]^{\b}\ar[r]^\pi&
X\ar[d]^{\a}\\ X'\ar[r]^{\e'}&Y'\ar[r]^{\pi'}& Y'\\ }
$$
such that $\pi\comp\e=\id_X$ and $\pi'\comp\e'=\id_{X'}$.
\end{rem}

Recall that an additive functor $F\colon\C\to\D$ is an {\em
  epimorphism of additive categories}, or simply an {\em epimorphism},
if $G\comp F=G'\comp F$ implies $G=G'$ for any pair $G,G'\colon\D\to
\E$ of additive functors.

\begin{lem}\label{le:mapsepi}
  Let $F\colon\C\to\D$ be an additive functor between additive
  categories having the following properties:
\begin{enumerate}
\item Every object in $\D$ belongs to the image of $F$.
\item For every map $\a\colon FX\to FY$ in $\D$, there are maps
$\a'\colon X'\to Y$ and $\pi\colon X'\to X$ in $\C$ such that
$F\a'=\a\comp F\pi$ and $F\pi$ is a split epimorphism.
\end{enumerate}
Then $F$ is an epimorphism.
\end{lem}
\begin{proof} 
  Let $G,G'\colon\D\to \E$ be a pair of additive functors satisfying
  $G\comp F=G'\comp F$.  The first condition implies that $G$ and $G'$
  coincide on objects, and the second condition implies that $G$ and
  $G'$ coincide on maps. Thus $G=G'$.
\end{proof}

Next we explain the notion of a flat functor.  A $\C^\op$-module $M$
is called {\em flat} if for every map $\b\colon Y\to Z$ in $\C$ and
every $y\in\Ker (M\b)$, there exists a map $\a\colon X\to Y$ in $\C$
and some $x\in MX$ such that $(M\a)x=y$ and $\b\comp\a=0$. We call an
additive functor $F\colon\C\to\D$ {\em flat} if the $\C^\op$-module
$\D(X,F-)$ is flat for every $X$ in $\D$.  We use the exact
structure of a module category in order to characterize flat functors.
Recall that $\Mod\C$ is an abelian category, and that a sequence $L\to
M\to N$ of $\C$-modules is exact if the sequence $LX\to MX\to NX$ is
exact for all $X$ in $\C$.  We record without proof a number of
equivalent conditions which justify our terminology.

\begin{lem}
  Let $F\colon\C\to\D$ be an additive functor between additive
  categories. Suppose $\mod\C$ is abelian. Then the following are
  equivalent.
\begin{enumerate}
\item $\D(X,F-)$ is a flat $\C^\op$-module for every $X$ in $\D$.
\item $F$ preserves weak kernels.
\item $F^\star \colon\mod\C\to\mod\D$ sends exact sequences to exact sequences.
\end{enumerate}
\end{lem}

\begin{lem}
  Let $F\colon\C\to\D$ be an additive functor between small additive
  categories. Then $F$ is flat if and only if
  $F^*\colon\Mod\C\to\Mod\D$ is an exact functor.
\end{lem}

Given an additive functor $F\colon\C\to\D$, we continue with a
criterion on $F$ such that $F^\star \colon\mod\C\to\mod\D$ is an exact
quotient functor.

\begin{lem}\label{le:flatepi}
  Let $F\colon\C\to\D$ be an additive functor between small additive
  categories. Suppose $\mod\C$ is abelian. If $F$
  is flat and $F_*\colon\Mod\D\to\Mod\C$ is fully faithful, then
  $\mod\D$ is abelian and $F^\star \colon\mod\C\to\mod\D$ is an exact
  quotient functor of abelian categories.
\end{lem}
\begin{proof}
  The functor $F^*\colon\Mod\C\to\Mod\D$ is exact because $F$ is flat,
  and it is a quotient functor because $F_*$ is fully faithful.  This
  follows from Lemma~\ref{le:abelian2} since $F_*$ is the right
  adjoint of $F^*$.  We conclude that $\mod\D$ is abelian and that the
  restriction $F^\star =F^*|_{\mod\C}$ to the category of finitely
  presented modules is an exact quotient functor, for instance by
  \cite[Theorem~2.6]{K1}.
\end{proof}

\begin{defn} 
Let $F\colon\C\to\D$ be an additive functor between additive
categories.  We call $F$ an {\em epimorphism up to direct factors}, if
there exists a  factorization $F=F_2\comp F_1$ such that 
\begin{enumerate}
\item $F_1$ is an epimorphism and bijective on objects, and
\item $F_2$ is fully faithful and every object in $\D$ is a direct
factor of some object in the image of $F_2$.
\end{enumerate}
\end{defn}

The following result summarizes our discussion and provides a
characterization of flat epimorphisms.

\begin{thm}\label{th:flatepi}
Let $F\colon\C\to\D$ be an additive functor between additive categories.
Suppose $\mod\C$ is abelian and $F$ is flat.
Then the following are equivalent.
\begin{enumerate}
\item The category $\mod\D$ is abelian and the exact functor $F^\star
  \colon\mod\C\to\mod\D$, sending $\C(-,X)$ to $\D(-,FX)$ for all $X$
  in $\C$, is a quotient functor of abelian categories.
\item Every object in $\D$ is a direct factor of some object in the
image of $F$. And for every map $\a\colon FX\to FY$ in $\D$, there are
maps $\a'\colon X'\to Y$ and $\pi\colon X'\to X$ in $\C$ such that
$F\a'=F\pi\comp\a$ and $F\pi$ is a split epimorphism.
\item $F$ is an epimorphism up to direct factors.
\end{enumerate}
\end{thm}
\begin{proof}
  (1) $\Rightarrow$ (2): Apply Lemma~\ref{le:maps}.

(2) $\Rightarrow$ (3): We define a factorization
$$\C\stackrel{F_1}\lto\D'\stackrel{F_2}\lto\D$$
as follows.  The
objects of $\D'$ are those of $\C$ and $F_1$ is the identity on
objects. Let $$\D'(X,Y)=\D(FX,FY)$$
for all $X,Y$ in $\C$, and let
$F_1\a=F\a$ for each map $\a$ in $\C$. The functor $F_2$ equals $F$ on
objects and is the identity on maps.  It follows from
Lemma~\ref{le:mapsepi} that $F_1$ is an epimorphism. The functor $F_2$
is fully faithful by construction, and $F_2$ is surjective up to
direct factors on objects by our assumption on $F$.

(3) $\Rightarrow$ (1): Assume that $F$ is an epimorphism up to direct
factors. We need to enlarge our universe so that $\C$ and $\D$ become
small categories. Note that this does not affect our assumption on
$F$, by Lemma~\ref{le:universe}. It follows from
Proposition~\ref{pr:epi4} that $F_*\colon \Mod\D\to\Mod\C$ is fully
faithful, and Lemma~\ref{le:flatepi} implies that $F^\star
\colon\mod\C\to\mod\D$ is a quotient functor.
\end{proof}

\section{Cohomological quotient functors}

In this section we introduce the concept of a cohomological quotient
functor between two triangulated categories. This concept generalizes
the classical notion of a quotient functor $\C\to\C/\B$ which Verdier
introduced for any triangulated subcategory $\B$ of $\C$; see
\cite{V}. 

\begin{defn}\label{de:coh}
Let $F\colon\C\to\D$ be an exact functor between triangulated
categories. We call $F$ a {\em cohomological quotient functor} if for
every cohomological functor $H\colon\C\to \A$ satisfying $\Ann
F\subseteq\Ann H$, there exists, up to a unique isomorphism, a unique
cohomological functor $H'\colon\D\to \A$ such that $H=H'\comp F$.
\end{defn}

Let us explain why a quotient funtor $\C\to\C/\B$ in the sense of
Verdier is a cohomological quotient functor. To this end we need the
following lemma.

\begin{lem}\label{le:ann}
Let $F\colon\C\to\D$ be an exact functor between triangulated
categories and suppose $F$ induces an equivalence $\C/\B\to\D$ for
some triangulated subcategory $\B$ of $\C$. Then $\Ann F$ is the ideal
of all maps in $\C$ which factor through some object in $\B$.
\end{lem}
\begin{proof}
The quotient $\C/\B$ is by definition the localization $\C[\Phi^{-1}]$
where $\Phi$ is the class of maps $X\to Y$ in $\C$ which fit into an
exact triangle $X\to Y\to Z\to\Si X$ with $Z$ in $\B$.  Now fix a map
$\psi\colon Y\to Z$ in $\Ann F$. The maps in $\C/\B$ are described via
a calculus of fractions. Thus $F\psi=0$ implies the existence of a map
$\p\colon X\to Y$ in $\Phi$ such that $\psi\comp\p=0$.  Complete
$\p$ to an exact triangle $X\to Y\to Z'\to\Si X$. Clearly, $\psi$ factors
through $Z'$ and $Z'$ belongs to $\B$. Thus $\Ann F$ is the ideal of
maps which factor through some object in $\B$.
\end{proof}

\begin{exm}\label{ex:verdier1}
A quotient functor $F\colon\C\to\C/\B$ is a cohomological quotient
functor. To see this, observe that a cohomological functor
$H\colon\C\to\A$ with $\Ker H$ containing $\B$ factors uniquely
through $F$ via some cohomological functor $H'\colon\C/\B\to\A$; see
\cite[Corollaire~II.2.2.11]{V}. Now use that
$$ \B\subseteq\Ker H\;\;\Longleftrightarrow\;\;\Ann F\subseteq\Ann H,$$
which follows from Lemma~\ref{le:ann}.
\end{exm}

It turns out that cohomological quotients are closely related to
quotients of additive and abelian categories.  The following result
makes this relation precise and provides a number of characterizations
for a functor to be a cohomological quotient functor.

\begin{thm}\label{th:cohquotient}
Let $F\colon\C\to\D$ be an exact functor between triangulated
categories.  Then the following are equivalent.
\begin{enumerate}
\item $F$ is a cohomological quotient functor.
\item The exact functor $F^\star \colon\mod\C\to\mod\D$, sending $\C(-,X)$ to
$\D(-,FX)$ for all $X$ in $\C$, is a quotient functor of abelian
categories.
\item Every object in $\D$ is a direct factor of some object in the
image of $F$. And for every map $\a\colon FX\to FY$ in $\D$, there are
maps $\a'\colon X'\to Y$ and $\pi\colon X'\to X$ in $\C$ such that
$F\a'=F\pi\comp\a$ and $F\pi$ is a split epimorphism.
\item $F$ is up to direct factors an epimorphism of additive categories. 
\end{enumerate}
\end{thm}
\begin{proof}
All we need to show is the equivalence of (1) and (2). The rest then
follows from Theorem~\ref{th:flatepi}.  Fix an abelian category $\A$
and consider the following commutative diagram
$$\xymatrix{
\HOM_\ex(\mod\D,\A)\ar[d]^{\HOM(H_\D,\A)}\ar[rr]^{\HOM(F^\star ,\A)}
&&\HOM_\ex(\mod\C,\A)\ar[d]^{\HOM(H_\C,\A)}\\
\HOM_\coh(\D,\A)\ar[rr]^{\HOM(F,\A)} &&\HOM_\coh(\C,\A) }$$ where the
vertical functors are equivalences.  Observe that $\HOM(H_\C,\A)$
identifies the exact functors $G\colon \mod\C\to\A$ satisfying $\Ker
F^\star \subseteq\Ker G$ with the cohomological functors $H\colon \C\to\A$
satisfying $\Ann F\subseteq\Ann H$.  This follows from the fact that
each $M$ in $\mod\C$ is of the form $M=\Im\C(-,\p)$ for some map $\p$
in $\C$. We conclude that the property of $F^\star $ to be an exact quotient
functor, is equivalent to the property of $F$ to be a cohomological
quotient functor.
\end{proof}

We complement the description of cohomological quotient functors by a
characterization of quotient functors in the sense of Verdier.

\begin{prop}\label{pr:verdier}
Let $F\colon\C\to\D$ be an exact functor between triangulated
categories.  Then the following are equivalent.
\begin{enumerate}
\item $F$ induces an equivalence $\C/{\Ker F}\to\D$.
\item Every object in $\D$ is isomorphic to some object in the
image of $F$. And for every map $\a\colon FX\to FY$ in $\D$, there are
maps $\a'\colon X'\to Y$ and $\pi\colon X'\to X$ in $\C$ such that
$F\a'=F\pi\comp\a$ and $F\pi$ is an isomorphism.
\end{enumerate}
\end{prop}
\begin{proof}
Let $\B=\Ker F$ and denote by $Q\colon\C\to \C/\B$ the quotient
functor, which is the identity on objects.  Given objects $X$ and $Y$
in $\C$, the maps $X\to Y$ in $\C/\B$ are fractions of the form
$$X\xrightarrow{(Q\pi)^{-1}} X'\xrightarrow{Q\a} Y$$ such that $F\pi$
is an isomorphism. 
This shows that (1) implies (2). To prove the converse, denote by
$G\colon\C/\B\to\D$ the functor which is induced by $F$. The
description of the maps in $\C/\B$ implies that $G$ is full. It
remains to show that $G$ is faithful. To this end choose a map
$\psi\colon Y\to Z$ such that $F\psi=0$. We complete $\psi$ to an
exact triangle
$$X\stackrel{\phi}\lto Y\stackrel{\psi}\lto Z\lto\Si X$$ and observe
that $F\p$ is a split epimorphism. Choose an inverse $\a\colon FY\to
FX$ and write it as $F\a'\comp (F\pi)^{-1}$, using (2). Thus
$Q(\phi\comp\a')$ is invertible, and $\psi\comp\phi\comp\a'=0$ implies
$Q\psi=0$ in $\C/\B$. We conclude that $G$ is faithful, and this
completes the proof.
\end{proof}

\section{Flat epimorphic quotients}

In this section we establish a triangulated structure for every
additive category which is a flat epimorphic quotient of some
triangulated category.

\begin{thm}\label{th:epi}
Let $\C$ be a triangulated category, and let $\D$ be an additive
category with split idempotents. Suppose $F\colon\C\to\D$ is a flat
epimorphism up to direct factors satisfying $\Si(\Ann F)=\Ann F$. Then
there exists a unique triangulated structure on $\D$ such that $F$ is
exact.  Moreover, a triangle $\Delta$ in $\D$ is exact if and only if there is an
exact triangle $\Gamma$ in $\C$ such that $\Delta$ is a direct factor
of $F\Gamma$.
\end{thm}

Note that an interesting application arises if one takes for $\D$ the
idempotent completion of $\C$. In this case, one obtains the main
result of \cite{BS}. 

The proof of Theorem~\ref{th:epi} is given in several steps and
requires some preparations. Assuming the suspension $\Si\colon\D\to\D$
is already defined, let us define the exact triangles in $\D$.  We
call a triangle $\Delta$ in $\D$ {\em exact}, if there exists an exact
triangle $\Gamma$ in $\C$ such that $\Delta$ is a direct factor of
$F\Gamma$, that is, there are triangle maps $\p\colon \Delta\to
F\Gamma$ and $\psi\colon F\Gamma\to\Delta$ such that
$\psi\comp\p=\id_\Delta$.

From now on assume that $F\colon\C\to\D$ is a flat epimorphism up to
direct factors, satisfying $\Si(\Ann F)=\Ann F$. We simplify our notation and
identify $\C$ with the image of the Yoneda functor $\C\to\mod\C$. The
same applies to the Yoneda functor $\D\to\mod\D$. Moreover, we
identify $F^\star=F$ and $\Si^\star=\Si$.  Note that $\mod\D$ is
abelian and that $F\colon\mod\C\to\mod\D$ is an exact quotient functor
by Theorem~\ref{th:flatepi}.  In particular, the maps in $\mod\D$ are
obtained from maps in $\mod\C$ via a calculus of fractions. Let us
construct the suspension for $\D$.

\begin{lem}
There is an equivalence $\Si'\colon\mod\D\to\mod\D$
making the following diagram commutative.
$$\xymatrix{ \mod\C\ar[d]^{\Si}\ar[r]^F&\mod\D\ar[d]^{\Si'}\\
\mod\C\ar[r]^F&\mod\D }
$$ The equivalence $\Si'$ is unique up to a unique isomorphism.
\end{lem}
\begin{proof}
Every object in $\mod\D$ is isomorphic to $FM$ for some $M$ in
$\mod\C$. And every map $\a\colon FM\to FN$ is a fraction,
that is, of the form
$$FM\xrightarrow{F\p} FN' \xrightarrow{(F\sigma)^{-1}} FN.$$ Now
define $\Si' (FM)=F(\Si M)$ and $\Si'\a=F(\Si\sigma)^{-1}\comp
F(\Si\p)$.
\end{proof}

We shall abuse notation and identify $\Si'=\Si$.  Now fix $M,N$ in
$\mod\D$. We may assume that $M=FM'$ and $N=FN'$.  We have a natural
map
$$\k_{M',N'}\colon\Hom_\C(M',N')\lto\Ext^3_\C(\Si M',N')$$ which is
induced from the triangulated structure on $\C$; see Appendix~\ref{ap:abel}.
This map induces a natural map
$$\k_{M,N}\colon\Hom_\D(M,N)\lto\Ext^3_\D(\Si M,N)$$ since every map
$FM'\to FN'$ is a fraction of maps in the image of $F$. Recall that
$\k_M=\k_{M,M}(\id_M)$. Let $$\Delta\colon X\stackrel{\a}\lto
Y\stackrel{\b}\lto Z\stackrel{\g}\lto\Si X$$ be a triangle in $\D$ and
put $M=\Ker\a$. We call $\Delta$ {\em pre-exact}, if 
$\g$ induces a map $Z\to\Si M$ such that the sequence
$$0\lto M\lto X\stackrel{\a}\lto Y\stackrel{\b}\lto Z\lto \Si M\lto
0$$ is exact in $\mod\D$ and represents $\k_M\in\Ext_\D^3(\Si M,M)$.
Note that every exact triangle in $\D$ is pre-exact, by
Proposition~\ref{pr:kappa}.

\begin{lem}\label{le:idpt}
Given a commutative diagram
\begin{equation}\label{eq:map1}
\xymatrix{
X\ar[d]^\p\ar[r]^\a&Y\ar[d]^\psi\ar[r]^\b&Z\ar[r]^\g&\Si X\ar[d]^{\Si\p}\\
X'\ar[r]^{\a'}&Y'\ar[r]^{\b'}&Z'\ar[r]^{\g'}&\Si X'\\
}
\end{equation}
in $\D$ such that both rows are pre-exact triangles, there exists a
map $\rho\colon Z\to Z'$ such that the completed diagram commutes.
Moreover, if $\p^2=\p$ and $\psi^2=\psi$, then there exists a choice
for $\rho$ such that $\rho^2=\rho$.
\end{lem}
\begin{proof}
Let $M=\Ker\a$ and $M'=\Ker\a'$. The pair $\p,\psi$ induces a map
$\mu\colon M\to M'$ and we obtain the following diagram in $\mod\D$.
\begin{equation}\label{eq:map2}
\xymatrix{ \k_M\colon
&0\ar[r]&\Omega^{-2}M\ar[d]^{\Omega^{-2}\mu}\ar[r]&Z\ar[r]&\Si
M\ar[r]\ar[d]^{\Si\mu}&0\\ \k_{M'}\colon
&0\ar[r]&\Omega^{-2}M'\ar[r]&Z'\ar[r]&\Si M'\ar[r]&0 }
\end{equation}
Here we use a dimension shift to represent $\k_M$ and $\k_{M'}$ by
short exact sequences.  The map $\k_{M,N}$ is natural in $M$ and $N$,
and therefore $\k_{M,\mu}(\k_M)=\k_{\mu,M'}(\k_{M'})$. This implies
the existence of a map $\rho\colon Z\to Z'$ making the diagram
(\ref{eq:map2}) commutative. Note that we can choose $\rho$ to be
idempotent if $\mu$ is idempotent.  It follows that the map $\rho$
completes the diagram (\ref{eq:map1}) to a map of triangles.
\end{proof}

\begin{lem}\label{le:complete}
Every map $X\to Y$ in $\D$ can be completed to an exact triangle
$X\to Y\to Z\to\Si X$.
\end{lem}
\begin{proof}
A map in $\D$ is a direct factor of some map in the image of $F$ by
Theorem~\ref{th:flatepi}; see also Remark~\ref{re:maps}. Thus we have
a commutative square
$$\xymatrix{
FX'\ar[r]^{F\a}\ar[d]^\p&FY'\ar[d]^\psi\\
FX'\ar[r]^{F\a}&FY'\\
}
$$ 
such that $\p$ and $\psi$ are idempotent and the map $X\to Y$ equals
the map $\Im\p\to\Im\psi$ induced by $F\a$. We complete $\a$ to an
exact triangle $\Delta$ in $\C$ and extend the pair $\p,\psi$ to an idempotent
triangle map $\e\colon F\Delta\to F\Delta$, which is possible by
Lemma~\ref{le:idpt}.  The image $\Im\e$ is an exact triangle in $\D$,
which completes the map $X\to Y$.
\end{proof}

We are now in the position to prove the octahedral axiom for $\D$.
Note that we have already established that $\D$ is a pre-triangulated
category. We say that a pair of composable maps $\a\colon X\to Y$ and
$\b\colon Y\to Z$ can be {\em completed to an octahedron} if
there exists a commutative diagram of the form
$$
\xymatrix{
X\ar@{=}[d]\ar[r]^\a&Y\ar[d]^\b\ar[r]&U\ar[r]\ar[d]&\Si X\ar@{=}[d]\\
X\ar[r]^{\b\comp\a}&Z\ar[r]\ar[d]^\g&V\ar[r]\ar[d]&\Si X\ar[d]^{\Si\a}\\
&W\ar@{=}[r]\ar[d]^\d&W\ar[d]\ar[r]^\d&\Si Y\\
&\Si Y\ar[r]&\Si U
}
$$
such that all triangles which occur are exact.

We shall use the following result due to Balmer and Schlichting.

\begin{lem}\label{le:BS} 
Let $\a\colon X\to Y$ and $\b\colon Y\to Z$ be maps in a
pre-triangulated category. Suppose there are objects $X',Y',Z'$ such
that
$$X\amalg X'
\xrightarrow{\left[\begin{smallmatrix}\a&0\\0&0\end{smallmatrix}\right]}
Y\amalg Y'\quad\textrm{and}\quad Y\amalg Y'
\xrightarrow{\left[\begin{smallmatrix}\b&0\\0&0\end{smallmatrix}\right]}
Z\amalg Z'$$ can be completed to an octahedron.  Then $\a$ and $\b$
can be completed to an octahedron.
\end{lem}
\begin{proof} See the proof of Theorem~1.12 in \cite{BS}.
\end{proof}

\begin{lem}\label{le:oct}
Every pair of composable maps in $\D$ can be completed to an octahedron.
\end{lem}
\begin{proof}
Fix two maps $\a\colon X\to Y$ and $\b\colon Y\to Z$ in $\D$.  We
proceed in two steps. First assume that $X=FA$, $Y=FB$, and $Z=FC$.
We use the description of the maps in $\D$ which is given in
Theorem~\ref{th:flatepi}.  We consider the map $\b\colon Y\to Z$ and
obtain new maps $\psi\colon B'\to C$ and $\pi\colon B'\to B$ in $\C$
such that $F\psi=\beta\comp F\pi$ and $F\pi$ is a split
epimorphism. We get a decomposition $FB'=Y\amalg Y'$ and an
automorphism $\e\colon Y\amalg Y'\to Y\amalg Y'$ such that
$F\psi\comp\e=\left[ \begin{smallmatrix}
\b&0\end{smallmatrix}\right]$. The same argument, applied to the
composite
$$X\stackrel{\a}\lto Y\stackrel{\left[
\begin{smallmatrix}\id_Y\\0\end{smallmatrix}\right]}\lto Y\amalg Y',$$ 
gives a map $\p\colon A'\to B'$ in $\C$, a decomposition $FA'=X\amalg
X'$, and an automorphism $\d\colon X\amalg X'\to X\amalg X'$ such that
$F\p\comp\d= \left[ \begin{smallmatrix}
\a&0\\0&0\end{smallmatrix}\right]$.  We know that the pair $\p,\psi$
in $\C$ can be completed to an octahedron.  Thus $F\p$ and $F\psi$ can
be completed to an octahedron in $\D$. It follows that $\left[
\begin{smallmatrix} \a&0\\0&0\end{smallmatrix}\right]$ and $\left[
\begin{smallmatrix} \b&0\end{smallmatrix}\right]$ can be completed to
an octahedron.  Using Lemma~\ref{le:BS}, we conclude that the pair
$\a,\b$ can be completed to an octahedron.

In the second step of the proof, we assume that the objects $X$, $Y$,
and $Z$ are arbitrary. Applying again the description of the maps in
$\D$, we find objects $X'$, $Y'$, and $Z'$ in $\D$ such that $X\amalg
X'$, $Y\amalg Y'$, and $Z\amalg Z'$ belong to the image of $F$. We
know from the first part of the proof that the maps
$$X\amalg X'
\xrightarrow{\left[\begin{smallmatrix}\a&0\\0&0\end{smallmatrix}\right]}
Y\amalg Y'\quad\textrm{and}\quad Y\amalg Y'
\xrightarrow{\left[\begin{smallmatrix}\b&0\\0&0\end{smallmatrix}\right]}
Z\amalg Z'$$ can be completed to an octahedron.  From this it follows
that $\a$ and $\b$ can be completed to an octahedron, using again
Lemma~\ref{le:BS}. This finishes the proof of the octahedral axiom for
$\D$.
\end{proof}

Let us complete the proof of Theorem~\ref{th:epi}.
\begin{proof}[Proof of Theorem~\ref{th:epi}]
We have constructed an equivalence $\Si\colon\D\to\D$, and the exact
triangles in $\D$ are defined as well.  We need to verify the axioms
(TR1) -- (TR4) from \cite{V}.  Let us concentrate on the properties of
$\D$, which are not immediately clear from our set-up.  In
Lemma~\ref{le:complete}, it is shown that every map in $\D$ can be
completed to an exact triangle.  In Lemma~\ref{le:idpt}, it is shown
that every partial map between exact triangles can be completed to a
full map. Finally, the octahedral axiom (TR4) is established in
Lemma~\ref{le:oct}.
\end{proof}

\section{A criterion for exactness}

Given an additive functor $\C\to\D$ between triangulated categories,
it is a natural question to ask when this functor is exact. We provide
a criterion in terms of the induced functor $\mod\C\to\mod\D$ and the
extension $\k_M$ in $\Ext^3(\Si^\star M,M)$ defined for each $M$ in
$\mod\C$; see Appendix~\ref{ap:abel}. There is an interesting
consequence. Given any factorization $F=F_2\comp F_1$ of an exact
functor, the functor $F_2$ is exact provided that $F_1$ is a
cohomological quotient functor.

\begin{prop}\label{pr:triangle}
Let $F\colon\C\to\D$ be an additive functor between triangulated
categories. Then $F$ is exact if and only if the following holds:
\begin{enumerate}
\item The right exact functor $F^\star \colon\mod\C\to\mod\D$, sending
$\C(-,X)$ to $\D(-,FX)$ for all $X$ in $\C$, is exact.
\item There is a natural isomorphism $\eta\colon
F\comp\Si_\C\to\Si_\D\comp F$.
\item $F^\star \k_M=\Ext^3_\D(\eta^\star_M,F^\star M)(\k_{F^\star M})$
for all $M$ in $\mod\C$.
\end{enumerate}
\end{prop}
Here, we denote by $\eta^\star$ the natural isomorphism $F^\star
\comp\Si^\star_\C\to\Si^\star_\D\comp F^\star$ which extends $\eta$,
that is, $\eta^\star_{\C(-,X)}=\D(-,\eta_X)$ for all $X$ in $\C$.

\begin{proof}
Suppose first that (1) -- (3) hold.  Let $$\Delta\colon X
\stackrel{\alpha}\longrightarrow Y\longrightarrow
Z\longrightarrow\Si_\C X$$ be an exact triangle in $\C$.  We need
to show that $F$ sends this triangle to an exact triangle in $\D$. To
this end complete the map $F\alpha$ to an exact triangle
$$ FX\stackrel{F\alpha}\longrightarrow FY\longrightarrow
Z'\longrightarrow\Si_\D(FX).$$ Now let $M=\Ker\C(-,\a)$.  We use
a dimension shift to represent the class $\k_M$ by a short exact
sequence corresponding to an element in $\Ext_\C^1(\Si^\star_\C
M,\Omega^{-2}M)$. Analogously, we represent $\k_{F^\star M}$ by a short
exact sequence. Next we use the exactness of $F^\star $ to obtain the
following diagram in $\mod\D$.
$$\xymatrix{F^\star \k_{M}\colon& 0\ar[r]&F^\star
(\Omega^{-2}M)\ar[r]\ar@{=}[d]&\D(-,FZ)\ar[r]& F^\star (\Si^\star_\C
M)\ar[r]\ar[d]^{\eta^\star_M}&0\\ \k_{F^\star M}\colon&
0\ar[r]&\Omega^{-2}(F^\star M)\ar[r]&\D(-,Z')\ar[r]&\Si^\star_\D
(F^\star M)\ar[r]&0\\ }$$ The diagram can be completed by a map
$\D(-,FZ)\to\D(-,Z')$ because $F^\star \k_M=\Ext^3_\D(\eta^\star_M,F^\star
M)(\k_{F^\star M})$. Let $\p\colon FZ\to Z'$ be the new map which is
an isomorphism, since $\eta^\star_M$ is an isomorphism.  We obtain the
following commutative diagram
$$\xymatrix{FX\ar[r]^{F\alpha}\ar@{=}[d]&FY\ar[r]\ar@{=}[d]&
FZ\ar[r]\ar[d]^\p&F(\Si_\C X)\ar[d]^{\eta_X}\\
FX\ar[r]^{F\alpha}&FY\ar[r]&Z'\ar[r]&\Si_\D(F X)}$$ and therefore
$F\Delta$ is an exact triangle.  Thus $F$ is an exact functor.  It
is not difficult to show that an exact functor $F$ satisfies (1) --
(3), and therefore the proof is complete.
\end{proof}

\begin{cor}\label{co:factorization}
Let $F\colon\C\to\D$ and $G\colon\D\to\E$ be additive functors between
triangulated categories. Suppose $F$ and $G\comp F$ are exact.
Suppose in addition that $F$ is a cohomological quotient functor. Then
$G$ is exact.
\end{cor}
\begin{proof}
We apply Proposition~\ref{pr:triangle}. First observe that $G^\star
\colon\mod\D\to\mod\E$ is exact because the composite $G^\star F^\star
=(GF)^\star$ is exact and $F^\star $ is an exact quotient functor, by
Theorem~\ref{th:cohquotient}.  Denote by $\eta^F\colon F^\star
\Si^\star_\C\to \Si^\star_\D F^\star $ and $\eta^{GF}\colon (G^\star F^\star
)\Si^\star_\C\to \Si^\star_\E (G^\star F^\star )$ the natural isomorphisms
which exists because $F$ and $GF$ are exact. In order to define
$\eta^G\colon G^\star \Si^\star_\D\to \Si^\star_\E G^\star $, we use again
the fact that $F^\star \colon\mod\C\to\mod\D$ is an exact quotient
functor. Thus every object in $\mod\D$ is isomorphic to $F^\star M$
for some $M$ in $\mod\C$. Moreover, any morphism $F^\star M\to F^\star
N$ is a fraction, that is, of the form
$$F^\star M\xrightarrow{F^\star \p} F^\star N' \xrightarrow{(F^\star
\sigma)^{-1}} F^\star N.$$ Now define $\eta^G_{F^\star M}$ as the
composite
$$\eta^G_{F^\star M}\colon (G^\star \Si^\star_\D)F^\star
M\xrightarrow{(G^\star \eta_M^F)^{-1}} (G^\star F^\star \Si^\star_\C)
M\xrightarrow{\eta_M^{GF}} (\Si^\star_\E G^\star )F^\star M.$$ The map
is natural, because $\eta^F$ and $\eta^{GF}$ are natural
transformations, and maps $F^\star M\to F^\star N$ come from maps in
$\mod\C$.  A straightforward calculation shows that $G^\star
\k_{N}=\Ext^3_\D(\eta^G_{N},G^\star N)(\k_{G^\star N})$ for all
$N=F^\star M$ in $\mod\D$.  Thus $G$ is exact by
Proposition~\ref{pr:triangle}.
\end{proof}

\section{Exact quotient functors}
The definition of a cohomological quotient functor between two
triangulated categories involves cohomological functors to an abelian
category.  It is natural to study the analogue where the cohomological
functors are replaced by exact functors to a triangulated category.

\begin{defn}\label{de:exactq}
Let $F\colon\C\to\D$ be an exact functor between triangulated
categories.  We call $F$ an {\em exact quotient functor} if for every
triangulated category $\E$ and every exact functor $G\colon\C\to\E$
satisfying $\Ann F\subseteq\Ann G$, there exists, up to a unique
isomorphism, a unique exact functor $G'\colon\D\to\E$ such that
$G=G'\comp F$.
\end{defn}

The motivating examples for this definition are the quotient functors
in the sense of Verdier.

\begin{exm}\label{ex:verdier2}
A quotient functor $F\colon\C\to\C/\B$ is an exact quotient
functor. To see this, observe that an exact functor $G\colon\C\to\D$
with $\Ker G$ containing $\B$ factors uniquely through $F$ via some
exact functor $G'\colon\C/\B\to\D$; see
\cite[Corollaire~II.2.2.11]{V}. Now use that
$$ \B\subseteq\Ker G\;\;\Longleftrightarrow\;\;\Ann F\subseteq\Ann G,$$
which follows from Lemma~\ref{le:ann}.
\end{exm}

We want to relate cohomological and exact quotient functors.

\begin{lem}\label{le:dominant}
Let $F\colon\C\to\D$ be a cohomological quotient functor, and denote
by $\D'$ the smallest full triangulated subcategory containing the
image of $F$. Then the restriction $F'\colon\C\to\D'$ of $F$ has the
following properties.
\begin{enumerate}
\item $F'$ is a cohomological quotient functor.
\item $F'$ is an exact quotient functor.
\end{enumerate}
\end{lem}
\begin{proof}
(1) Use the characterization of cohomological quotient functors in
Theorem~\ref{th:cohquotient}.

(2) For simplicity we assume $\D'=\D$. Let $G\colon\C\to\E$ be an
exact functor satisfying $\Ann F\subseteq\Ann G$. Then the composite
$H_{\E}\comp G$ with the Yoneda functor factors through $F$ because
$F$ is a cohomological quotient functor. We have the following
sequence of inclusions
$$\E\subseteq\bar\E\subseteq\mod\E$$ where $\bar\E$ denotes the
idempotent completion of $\E$.  We obtain a functor $\D\to\mod\E$ and
its image lies in $\bar\E$, since every object in $\D$ is a direct
factor of some object in the image of $F$. Thus we have a functor
$G'\colon\D\to\bar\E$ which is exact by
Corollary~\ref{co:factorization}. Our additional assumption on $F$
implies that $\Im G'\subseteq \E$. We conclude that $G$ factors
through $F$ via an exact functor $\D\to\E$.
\end{proof}

The following example has been suggested by B. Keller. It shows that
there are exact quotient functors which are not cohomological quotient
functors.

\begin{exm}
  Let $A$ be the algebra of upper $2\times 2$ matrices over a field
  $k$, and let $B=k\times k$.  We consider the bounded derived
  categories $\C=\mathbf D^b(\mod A)$ and $\D=\mathbf D^b(\mod B)$.
  Restriction along the algebra homomorphism $f\colon B\to A$,
  $(x,y)\mapsto\left[\begin{smallmatrix}x&0\\
      0&y\end{smallmatrix}\right]$, induces an exact functor $F\colon
  \C\to\D$ which is an exact quotient functor but not a cohomological
  quotient functor. In fact, $f$ has a left inverse $A\to B$,
  $\left[\begin{smallmatrix}x&z\\ 0&y\end{smallmatrix}\right]\mapsto
  (x,y)$, which induces a right inverse $G\colon \D\to\C$ for $F$.
  Thus every exact functor $F'\colon\C\to\E$ satisfying $\Ann
  F\subseteq \Ann F'$ factors uniquely through $F$, by
  Lemma~\ref{le:rightinverse} below.  However, the exact functor
  $F^\star\colon \mod\C\to\mod\D$ extending $F$ does not induce an
  equivalence $\mod\C/{\Ker F^\star}\to\mod\D$.

Let us describe $\mod\C$. To this end denote by $0\to X_1\to X_3\to
X_2\to 0$ the unique non-split exact sequence in $\mod A$ involving
the simple $A$-modules $X_1$ and $X_2$. This sequence induces an exact
triangle $X_1\to X_3\to X_2\to \Si X_1$ in $\C$. Note that each
indecomposable object in $\C$ is determined by its cohomology and is
therefore of the form $\Si^n X_i$ for some $n\in\mathbb Z$ and some
$i\in\{1,2,3\}$. Thus the indecomposable objects in $\mod\C$ are
precisely the objects of the form
$$\C(-,\Si^nX_i)/{\rad^j\C(-,\Si^nX_i)}\quad \textrm{with $n\in\mathbb
Z$, $i\in\{1,2,3\}$, $j\in\{0,1\}$},$$ where $\rad^0M=M$ and $\rad^1M$
is the intersection of all maximal subobjects of $M$.  The restriction
functor $\mod A\to\mod B$ sends $0\to X_1\to X_3\to X_2\to 0$ to a
split exact sequence. Thus $F$ kills the map $X_2\to\Si X_1$ in $\C$,
and we have that $F^\star M=0$ for some indecomposable $M$ in $\mod\C$ if
and only if $M\cong \C(-,\Si^nX_2)/{\rad\C(-,\Si^nX_2)}$ for some
$n\in\mathbb Z$.  It follows that the canonical functor
$$\coprod_{n\in\mathbb Z}\mod A\lto\mod\C/{\Ker
F^\star},\quad (M_n)_{n\in\mathbb Z}\longmapsto\coprod_{n\in\mathbb
Z}\C(-,\Si^nM_n)$$ is an equivalence.

We have seen that $\mod\C/{\Ker F^\star}$ is not a semi-simple category,
whereas in $\mod\D$ every object is semi-simple. More specifically,
the cohomological functor
$$H\colon\C\lto\mod\C\lto\mod\C/{\Ker F^\star}$$ does not factor through
$F$ via some cohomological functor $\D\to\mod\C/{\Ker F^\star}$, even
though $\Ann H=\Ann F$.  We conclude that $F$ is not a cohomological
quotient functor.
\end{exm}

\begin{lem}\label{le:rightinverse}
  Let $F\colon\C\to\D$ be an additive functor between additive
  categories which admits a right inverse $G\colon\D\to\C$, that is,
  $F\comp G=\Id_\D$. Suppose $F'\colon\C\to\E$ is an additive functor
  satisfying
\begin{enumerate}
\item $\Ann F\subseteq\Ann F'$, and
\item for all $X,Y$ in $\C$, $FX=FY$ implies $F'X=F'Y$. 
\end{enumerate}
Then $F'$ factors uniquely through $F$.
\end{lem}
\begin{proof}
  We have $F\comp(G\comp F-\Id_\C)=0$ and therefore $F'\comp(G\comp
  F-\Id_\C)=0$.  Thus $F'=(F'\comp G)\comp F$.
\end{proof}

\section{Exact ideals}

Given a cohomological quotient functor $F\colon\C\to\D$, the ideal
$\Ann F$ is an important invariant. In this section we investigate the
collection of all ideals which are of this form.

\begin{defn}\label{de:exact} 
Let $\C$ be a triangulated category.  An ideal $\mathfrak I$ of $\C$
is called {\em exact} if there exists a cohomological quotient functor
$F\colon\C\to\D$ such that $\mathfrak I=\Ann F$.
\end{defn}

The exact ideals are partially ordered by inclusion and we shall
investigate the structure of this poset. Recall that an ideal
$\mathfrak I$ in a triangulated category $\C$ is {\em cohomological},
if there exists a cohomological functor $F\colon\C\to\A$ such that
$\mathfrak I=\Ann F$.

\begin{thm}\label{th:ideallattice}
Let $\C$ be a small triangulated category. Then the exact ideals in
$\C$ form a complete lattice, that is, given a family $(\mathfrak
I_i)_{i\in\La}$ of exact ideals, the supremum $\sup \mathfrak I_i$ and
the infimum $\inf \mathfrak I_i$ exist. Moreover, the supremum
coincides with the supremum in the lattice of cohomological ideals.
\end{thm}

Our strategy for the proof is to use a bijection between the
cohomological ideals of $\C$ and the Serre subcategories of
$\mod\C$. We proceed in several steps and start with a few
definitions.  Given an ideal $\mathfrak I$ of $\C$, we define
$$\Im\mathfrak I=\{M\in\mod\C\mid \mbox{$M\cong\Im\C(-,\p)$ for some
  $\p\in\mathfrak I$}\}.$$
The next definition is taken from \cite{Be}.
\begin{defn}\label{de:saturated} 
Let $\C$ be a triangulated category.  An ideal $\mathfrak I$ of $\C$
is called {\em saturated} if for every exact triangle
$X'\stackrel{\a}\to X\stackrel{\b}\to X''\to\Si X'$ and every map
$\p\colon X\to Y$ in $\S_c$, we have that $\p\comp\a,\b\in\mathfrak
I$ implies $\p\in\mathfrak I$.
\end{defn}  
The following characterization combines \cite[Lemma~3.2]{K} and
\cite[Theorem~3.1]{Be}.

\begin{lem}\label{le:serre} 
Let $\C$ be a triangulated category. Then the following are equivalent
for an ideal $\mathfrak I$ of $\C$.
\begin{enumerate}
\item $\mathfrak I$ is cohomological.
\item $\mathfrak I$ is saturated.
\item $\Im \mathfrak I$ is a Serre subcategory of $\mod\C$.
\end{enumerate}
Moreover, the map $\mathfrak J\mapsto\Im \mathfrak J$ induces a
bijection between the cohomological ideals of $\C$ and the Serre
subcategories of $\mod\C$.
\end{lem}
\begin{proof} 
(1) $\Rightarrow$ (2): Let $\mathfrak I=\Ann F$ for some cohomological
functor $F\colon \C\to\A$.  Fix an exact triangle $X'\stackrel{\a}\to
X\stackrel{\b}\to X''\to\Si X'$ and a map $\p\colon X\to Y$ in
$\C$. Suppose $\p\comp\a,\b\in\mathfrak I$. Then $F\a$ is an
epimorphism, and therefore $F\p\comp F\a=0$ implies $F\p=0$.  Thus
$\p\in\mathfrak I$.

(2) $\Rightarrow$ (3): Let $0\to F'\to F\to F''\to 0$ be an exact
sequence in $\mod\C$. Using that $\mathfrak I$ is an ideal, it is
clear that $F\in\Im\mathfrak I$ implies $F',F''\in\Im\mathfrak I$.
Now suppose that $F',F''\in\Im\mathfrak I$. Using that $\mod\C$ is a
Frobenius category, we find maps $\p\colon X\to Y$ and $\a\colon
X'\to X$ such that $F=\Im\C(-,\p)$ and $F'=\Im\C(-,\p\comp\a)$. Now
form exact triangles 
$$W\stackrel{\chi}\lto X\stackrel{\p}\lto Y\lto\Si W\quad\textrm{and}\quad
X'\amalg
W\stackrel{\left[\begin{smallmatrix}\a&\chi\end{smallmatrix}\right]}\lto
X\stackrel{\b}\lto X''\lto\Si (X'\amalg W)$$
in $\C$, and observe that
$F''=\Im\C(-,\b)$.  We have
$\p\comp\left[\begin{smallmatrix}\a&\chi\end{smallmatrix}\right]$ and
$\b$ in $\mathfrak I$.  Thus $\p\in\mathfrak I$ since $\mathfrak I$
is saturated. It follows that $F=\Im\C(-,\p)$ belongs to
$\Im\mathfrak I$.

(3) $\Rightarrow$ (1): Let $F$ be the composite of the Yoneda functor
$\C\to\mod\C$ with the quotient functor $\mod\C\to\mod\C/{\Im\mathfrak
  I}$. This functor is cohomological and we have $\mathfrak I=\Ann F$.
\end{proof}

We need some more terminology. Fix an abelian category $\A$. A Serre
subcategory $\B$ of $\A$ is called {\em localizing} if the quotient
functor $\A\to\A/\B$ has a right adjoint. If $\A$ is a Grothendieck
category, then $\B$ is localizing if and only if $\B$ is closed under
taking coproducts \cite[Proposition~III.8]{G}. We denote for any
subcategory $\B$ by $\li\B$ the full subcategory of filtered colimits
$\li X_i$ in $\A$ such that $X_i$ belongs to $\B$ for all $i$.

Now let $\C$ be a small additive category and suppose $\mod\C$ is
abelian.  Given a Serre subcategory $\S$ of $\mod\C$, then $\li\S$ is
a localizing subcategory of $\Mod\C$; see \cite[Theorem~2.8]{K1}. 
This has the following consequence which we record for later
reference.

\begin{lem}\label{le:localizing}
Let $\C$ be a small triangulated category and $\mathfrak I$ be a
cohomological ideal of $\C$. Then $\li\Im\mathfrak I$ is a localizing
subcategory of $\Mod\C$.
\end{lem}
\begin{proof} Use Lemma~\ref{le:serre}.
\end{proof}

We call a Serre subcategory $\S$ of $\mod\C$ {\em perfect} if the right
adjoint of the quotient functor $\Mod\C\to\Mod\C/{\li\S}$ is an exact
functor. We have a correspondence between perfect Serre subcategories of
$\mod\C$ and flat epimorphisms starting in $\C$. To make this precise,
we call a pair $F_1\colon\C\to\D_1$ and $F_2\colon\C\to\D_2$ of flat
epimorphisms {\em equivalent} if $\Ker F_1^\star=\Ker F_2^\star$.

\begin{lem}\label{le:flatSerre}
  Let $\C$ be a small additive category and suppose $\mod\C$ is
  abelian.  Then the map
  $$(F\colon\C\to\D)\;\;\longmapsto\;\;\Ker F^\star$$
  induces
  a bijection between the equivalence classes of flat epimorphisms
  starting in $\C$, and the perfect Serre subcategories of $\mod\C$.
\end{lem}
\begin{proof}
  We construct the inverse map as follows. Let $\S$ be a perfect Serre
  subcategory of $\mod\C$ and consider the quotient functor
  $Q\colon\Mod\C\to\Mod\C/{\li\S}$. Observe that $Q$ preserves
  projectivity since the right adjoint of $Q$ is exact. Now define
  $\D$ to be the full subcategory formed by the objects $Q\C(-,X)$
  with $X$ in $\C$, and let $F\colon\C\to\D$ be the functor which sends $X$ to
  $Q\C(-,X)$. It follows that $F^*\colon\Mod\C\to\Mod\D$
  induces an equivalence $\Mod\C/{\li\S}\to\Mod\D$. Thus $F$ is a flat
  epimorphism satisfying $\S=\Ker F^\star$.
\end{proof}

\begin{lem}\label{le:serrelattice}
Let $\C$ be a small additive category and suppose $\mod\C$ is abelian.
If $(\S_i)_{i\in\La}$ is a family of perfect Serre subcategories of
$\mod\C$, then the smallest Serre subcategory of $\mod\C$ containing
all $\S_i$ is perfect.
\end{lem}
\begin{proof}
  For each $i\in\La$, let $\M_i$ be the full subcategory of
  $\C$-modules $M$ satisfying
  $\Hom_\C(\li\S_i,M)=0=\Ext^1_\C(\li\S_i,M)$.  Note that the right
  adjoint adjoint of the quotient functor
  $Q\colon\Mod\C\to\Mod\C/{\li\S_i}$ identifies $\Mod\C/{\li\S_i}$
  with $\M_i$; see Lemma~\ref{le:abelian2}. Let $\S=\sup\S_i$. Then
  the full subcategory $\li\S$ is the smallest localizing subcategory
  of $\Mod\C$ containing all $\S_i$. Let $\M$ be the full subcategory
  of $\C$-modules $M$ satisfying
  $\Hom_\C(\li\S,M)=0=\Ext^1_\C(\li\S,M)$.  We claim that
  $\M=\bigcap_i\M_i$. To see this, let $\I_i$ be the full subcategory
  of injective objects in $\M_i$. Note that a $\C$-module $M$ belongs
  to $\li\S_i$ if and only if $\Hom_\C(M,\I_i)=0$, and $M$ belongs to
  $\M_i$ if and only if the modules $I_0,I_1$ in a minimal injective
  resolution $0\to M\to I_0\to I_1$ belong to $\I_i$.  Let
  $\I=\bigcap_i\I_i$. Then we have that a $\C$-module $M$ belongs to
  $\li\S$ if and only if $\Hom_\C(M,\I)=0$, and $M$ belongs to
  $\bigcap_i\M_i$ if and only if the modules $I_0,I_1$ in a minimal
  injective resolution $0\to M\to I_0\to I_1$ belong to $\I$. This
  proves $\M=\bigcap_i\M_i$.  It follows that the inclusion
  $\M\to\Mod\C$ is exact because each inclusion $\M_i\to\Mod\C$ is
  exact. Thus $\S$ is a perfect Serre subcateory.
\end{proof}

\begin{prop}\label{pr:exactideal}
Let $\C$ be a small triangulated category, and let $\mathfrak I$ be an
ideal in $\C$ satisfying $\Si\mathfrak I=\mathfrak I$.  Then the
following are equivalent.
\begin{enumerate}
\item $\mathfrak I$ is an exact ideal. 
\item $\Im\mathfrak I$ is a perfect Serre subcategory of $\mod\C$.
\item There exists a flat epimorphism $F\colon\C\to\D$ such that $\Ann F=\mathfrak I$.
\end{enumerate}
\end{prop}
\begin{proof}
  (1) $\Rightarrow$ (2): Suppose $\mathfrak I$ is an exact ideal, that
  is, there is a cohomological quotient functor $F\colon\C\to\D$ such
  that $\mathfrak I=\Ann F$. Then $F$ is a flat epimorphism by
  Theorem~\ref{th:cohquotient}.  Now observe that $\Im\mathfrak I=\Ker
  F^\star$.  Thus $\Im\mathfrak I$ is perfect by
  Lemma~\ref{le:flatSerre}.
  
  (2) $\Rightarrow$ (3): Apply again Lemma~\ref{le:flatSerre} to
  obtain a flat epimorphisms $F\colon\C\to\D$ with $\Ann F=\mathfrak
  I$.
  
  (3) $\Rightarrow$ (1): We may assume that idempotents in $\D$ split.
  It follows from Theorem~\ref{th:epi} that $\D$ is a triangulated
  category and that $F$ is an exact functor. Moreover,
  Theorem~\ref{th:cohquotient} implies that $F$ is a cohomological
  quotient functor.  Thus $\mathfrak I$ is an exact ideal.
\end{proof}

We collect our findings to obtain the proof of the theorem from the
beginning of this section.

\begin{proof}[Proof of Theorem~\ref{th:ideallattice}.]
Let $(\mathfrak I_i)_{i\in\La}$ be a family of exact ideals in $\C$
and consider the corresponding Serre subcategories $\S_i=\Im\mathfrak
I_i$ of $\mod\C$ which are perfect by Proposition~\ref{pr:exactideal}.
It follows from Lemma~\ref{le:serrelattice} that $\S=\sup\S_i$ is
perfect.  There is a cohomological ideal $\mathfrak I$ in $\C$
satisfying $\S=\Im\mathfrak I$ and we have $\mathfrak I=\sup\mathfrak
I_i$ in the lattice of cohomological ideals by Lemma~\ref{le:serre}.
Applying again Proposition~\ref{pr:exactideal}, we see that the ideal
$\mathfrak I$ is exact. This completes the proof.
\end{proof}

\section{Factorizations}

Let $\C$ be a triangulated category. If $\B$ is a full triangulated
subcategory of $\C$, then the quotient functor $Q\colon\C\to\C/\B$ in
the sense of Verdier is a cohomological and exact quotient functor (in
the sense of Definition~\ref{de:coh} and
Definition~\ref{de:exactq}). This fact motivates the following
definition.

\begin{defn} 
Let $F\colon\C\to\D$ be an exact functor between triangulated
categories. We call $F$ a {\em CE-quotient functor} if $F$ is a
cohomological quotient functor and an exact quotient
functor.\begin{footnote} {The terminology refers to the properties
`cohomological' and `exact'.  In addition, we wish to honour Cartan
and Eilenberg.}
\end{footnote}
\end{defn}

In this section we study the collection of all CE-quotient functors
starting in a fixed triangulated category.  Given a pair
$F_1\colon\C\to\D_1$ and $F_2\colon\C\to\D_2$ of CE-quotient functors, we
define
\begin{align*}
F_1\sim F_2 & \;\;\Longleftrightarrow\;\; \textrm{there exists an
equivalence $G\colon\D_1\to\D_2$ such that $F_2=G\comp F_1$,}\\
F_1\geq F_2 & \;\;\Longleftrightarrow\;\; \textrm{there exists an
exact functor $G\colon\D_1\to\D_2$ such that $F_2=G\comp F_1$.}
\end{align*}
We obtain a partial ordering on the equivalence classes of
CE-quotient functors, which may be rephrased as follows.
\begin{align*}
F_1\sim F_2 & \;\;\Longleftrightarrow\;\; \Ann F_1=\Ann F_2, \\
F_1\geq F_2 & \;\;\Longleftrightarrow\;\; \Ann F_1\subseteq\Ann F_2.
\end{align*}
The ideals of the form $\Ann F$  arising from CE-quotient functors
$F\colon\C\to\D$ form a complete lattice. This has been
established in Theorem~\ref{th:ideallattice}, and we obtain the
following immediate consequence.

\begin{thm}\label{th:quotlattice}
Let $\C$ be a small triangulated category. 
\begin{enumerate}
\item The equivalence classes of CE-quotient functors starting in
$\C$ form a complete lattice.
\item The assignment $F\mapsto\Ann F$ induces an anti-isomorphism
between the lattice of CE-quotient functors starting in $\C$ and the
lattice of exact ideals of $\C$.
\item Given a family $(F_i)_{i\in\La}$ of CE-quotient functors, we
have
$$\Ann(\inf_{i\in\La}F_i)=\sup_{i\in\La}(\Ann F_i)$$
where the supremum is taken in the lattice of cohomological ideals of $\C$.
\end{enumerate}
\end{thm}
\begin{proof}
The ideals of the form $\Ann F$ for some CE-quotient functor
$F\colon\C\to\D$ are precisely the exact ideals of $\C$. This follows
from Lemma~\ref{le:dominant}. Now apply Theorem~\ref{th:ideallattice}.
\end{proof}

The completeness of the CE-quotient functor lattice yields a canonical
factorization for every exact functor between two triangulated
categories.

\begin{cor}\label{co:factor}
Let $\C$ be a small triangulated category. Then every exact functor
$F\colon\C\to\D$ to a triangulated category $\D$ has a factorization
$$\xymatrix@=1.5em{\C\ar[rr]^Q&&\C'\ar[rr]^{F'} &&\D}$$
having the following properties:
\begin{enumerate}
\item $Q$ is a CE-quotient functor and $F'$ is exact.
\item Given a factorization
$$\xymatrix@=1.5em{\C\ar[rr]^{Q'}&&\C''\ar[rr]^{F''} &&\D}$$ of $F$ such that
$Q'$ is a CE-quotient functor and $F''$ is exact,
there exists, up to a unique isomorphism, a unique exact functor
$G\colon \C''\to\C'$
$$\xymatrix@=1.5em{
&&\C''\ar[rrd]^{F''}\ar[dd]^G \\
\C\ar[rru]^{Q'}\ar[rrd]^{Q}&&&&\D\\
&&\C'\ar[rru]^{F'}
}$$
such that $Q=G\comp Q'$ and $F''\cong F'\comp G$.
\end{enumerate}
\end{cor}
\begin{proof}
We obtain the CE-quotient functor
$Q\colon\C\to\C'$ by taking the infimum over all CE-quotient functors
$Q'\colon \C\to\C''$ admitting a factorization
$$\xymatrix@=1.5em{\C\ar[rr]^{Q'}&&\C''\ar[rr]^{F''} &&\D.}$$ Note
that $F$ factors through $Q$ because $\Ann Q\subseteq \Ann F$.  This
follows from the fact that $\Ann F$ is cohomological and $\Ann
Q'\subseteq \Ann F$ for all $Q'$.
\end{proof}

\section{Compactly generated triangulated categories and Brown representability}

We recall the definition of a compactly generated triangulated
category, and we review a variant of Brown's Representability Theorem
which will be needed later on.

Let $\S$ be a triangulated category and suppose that arbitrary
coproducts exist in $\S$.  An object $X$ in $\S$ is called {\em
compact} if for every family $(Y_i)_{i\in I}$ in $\S$, the canonical
map $\coprod_i\S(X,Y_i)\to\S(X,\coprod_i Y_i)$ is an isomorphism.  We
denote by $\S_c$ the full subcategory of compact objects in $\S$ and
observe that $\S_c$ is a triangulated subcategory of $\S$.  Following
\cite{N1}, the category $\S$ is called {\em compactly generated}
provided that the isomorphism classes of objects in $\S_c$ form a set,
and $\S(C,X)=0$ for all $C$ in $\S_c$ implies $X=0$ for every object
$X$ in $\S$.

A basic tool for studying  a compactly generated triangulated
category $\S$ is the cohomological functor
$$H_\S\colon\S\longrightarrow\Mod\S_c,\quad X\mapsto
H_X=\S(-,X)|_{\S_c}$$ which we call {\em restricted Yoneda functor}.
Our notation does not distinguish between the Yoneda functor
$H_\S\colon\S\to\mod\S$ and the restricted Yoneda functor. However,
the meaning of $H_\S$ and $H_X$ for some $X$ in $\S$ will be clear
from the context.


Next we recall from \cite{K2} a variant of Brown's Representability
Theorem \cite{Br}; see also \cite{Ke1,N1,Fr}.  Let $\S$ be a
triangulated category with arbitrary products.  An object $U$ in $\S$
is called a {\em perfect cogenerator} if $\S(X,U)=0$ implies $X=0$ for
every object $X$ in $\S$, and for every countable family of maps
$X_i\to Y_i$ in $\S$, the induced map
$$\S(\prod_iY_i,U)\longrightarrow\S(\prod_iX_i,U)$$ is surjective
provided that the map $\S(Y_i,U)\rightarrow\S(X_i,U)$ is surjective
for all $i$.

\begin{prop}[Brown representability]\label{pr:brown} 
Let $\S$ be a triangulated category with arbitrary products and a
perfect cogenerator $U$.
\begin{enumerate}
\item A functor $H\colon\S\to\Ab$ is cohomological and preserves all products 
if and only if $H\cong\S(X,-)$ for some object $X$ in $\S$.
\item $\S$ coincides with its smallest full triangulated subcategory
which contains $U$ and is closed under taking all products.
\end{enumerate}
\end{prop}
\begin{proof} See Theorem~A in \cite{K2}.
\end{proof}

There is an immediate consequence which we shall use.

\begin{cor}\label{co:brown} 
Let $\S$ be a triangulated category with arbitrary products and a
perfect cogenerator $U$. An exact functor $\S\to\T$ between triangulated
categories preserves all products if and only if it has a
left adjoint.
\end{cor}
\begin{proof} 
The left adjoint of a functor $F\colon\S\to\T$ sends an object $X$ in
$\T$ to the object in $\S$ representing $\T(X,F-)$.
\end{proof}

Take as an example for $\S$ a compactly generated triangulated
category. Then
$$U=\coprod_{C\in\S_c}C$$ is a perfect cogenerator for
$\S^\op$.  If $I$ is an injective cogenerator for $\Mod\S_c$, then the
object $V$ satisfying
$$\Hom_{\S_c}(H_\S-,I)\cong\S(-,V)$$ is a perfect cogenerator for
$\S$. Note that $V$ exists because $\S^\op$ is perfectly cogenerated.

\section{Smashing localizations}

We establish for any compactly generated triangulated category $\S$ a
bijective correspondence
between the smashing localizations of $\S$ and the cohomological
quotients of $\S_c$.

This result is divided into two parts. In this section we show
that any smashing localization induces a cohomological quotient. Let
us recall the relevant definitions.

An exact functor $F\colon\S\to\T$ between triangulated categories is a
{\em localization functor} if it has a right adjoint $G$ such that
$F\comp G\cong\Id_\T$. Note that the condition $F\comp G\cong\Id_\T$
is equivalent to the fact that $F$ induces an equivalence $\S/{\Ker
F}\to \T$, where $\S/{\Ker F}$ denotes quotient in the sense of
Verdier \cite{V}. It is often useful to identify a localization
functor $F\colon\S\to\T$ with the idempotent functor $L\colon\S\to\S$
defined by $L=G\comp F$.  The {\em $L$-acyclic} objects are those in
$\Ker F$ and the {\em $L$-local} object are those which are isomorphic
to some object in the image of $G$.  The localization $F$ is called
{\em smashing} if $G$ preserves all coproducts which exist in
$\T$.\begin{footnote}{If $\S$ carries a smash product $\wedge\colon
\S\times\S\to\S$ with unit $S$, then $LX=X\wedge LS$ provided $F$ is
smashing.}
\end{footnote}

\begin{thm}\label{th:smash1}
Let $\S$ be a compactly generated triangulated category and
$F\colon\S\to \T$ be an exact functor between triangulated
categories. Then $F$ is a smashing localization if and only if the
following holds:
\begin{enumerate}
\item $\T$ is a compactly generated triangulated category.
\item $F$ preserves coproducts.
\item $F$ induces a functor $F_c\colon \S_c\to\T_c$ which is a
cohomological quotient functor.
\end{enumerate}
\end{thm}

We need a few preparations before we can give the proof of this
result.

\begin{lem}\label{le:adj}
Let $\S$ be a compactly generated triangulated category and
$F\colon\S\to \T$ be an exact functor between triangulated
categories. Suppose $F$ preserves coproducts.  Then the right adjoint
of $F$ preserves coproducts if and only if $F$ preserves compactness.
\end{lem}
\begin{proof}
Combine the definition of compactness and the adjointness isomorphism;
see \cite[Theorem~5.1]{N1}.
\end{proof}

\begin{lem}\label{le:mod} 
  Let $F\colon\S\to\T$ be an exact functor between compactly generated
  triangulated categories. Suppose $F$ has a right adjoint
  $G\colon\T\to\S$ which preserves coproducts. Then the following
  diagram commutes.
$$\xymatrix{
\S\ar[d]^{H_\S}\ar[r]^F&\T\ar[d]^{H_\T}\ar[r]^G&\S\ar[d]^{H_\S}\\
\Mod\S_c\ar[r]^{(F_c)^*}&\Mod\T_c\ar[r]^{(F_c)_*}&\Mod\S_c\\
}
$$
\end{lem}
\begin{proof}
See Proposition~2.6 in \cite{K}.
\end{proof}

We are now in the position that we can prove the main result of this section.

\begin{proof}[Proof of Theorem~\ref{th:smash1}.]
Suppose first that $F$ is a smashing localization. Thus $F$ has a
right adjoint $G$ which preserves coproducts. This implies that $F$
induces a functor $F_c\colon \S_c\to\T_c$, by Lemma~\ref{le:adj}.
Using the adjointness formula $\T(FX,Y)\cong\S(X,GY)$, one sees that
$\T$ is generated by the image of $F_c$. Thus $\T$ is compactly
generated.  It remains to show that $F_c$ is a cohomological quotient
functor. To this end denote by $\M$ the class of $\T_c$-modules $M$ such
that the natural map $((F_c)^*\comp(F_c)_*)M\to M$ is an isomorphism.
Observe that $(F_c)^*\comp(F_c)_*$ composed with the Yoneda
embedding $\T_c\to\Mod\T_c$ equals the composite $H_\T \comp F\comp
G|_{\T_c}$, by Lemma~\ref{le:mod}.  Our assumption implies $F\comp
G\cong\Id_\T$, and therefore $\M$ contains all representable functors.
The composite $(F_c)^*\comp(F_c)_*$ preserves all colimits and therefore
$\M$ is closed under taking colimits. We conclude that
$$(F_c)^*\comp(F_c)_*\cong\Id_{\Mod\T_c}$$
since every module is a
colimit of representable functors. Thus $F_c$ is up to direct factors an
epimorphism by Proposition~\ref{pr:epi4}, and therefore a
cohomological quotient functor by Theorem~\ref{th:cohquotient}.

Now suppose that $F$ satisfies (1) -- (3). An application of Brown's
Representability Theorem shows that $F$ has a right adjoint $G$, since
$F$ preserves coproducts.  Moreover, $G$ preserves coproducts by
Lemma~\ref{le:adj}, since $F$ preserves compactness. It remains to
show that $F\comp G\cong\Id_\T$.  To this end denote by $\T'$ the class
of objects $X$ in $\T$ such that the natural map $(F\comp G)X\to X$ is
an isomorphism.  Our assumption on $F_c$ implies
$$(F_c)^*\comp(F_c)_*\cong\Id_{\Mod\T_c}.$$ Using again
Lemma~\ref{le:mod}, we see that $\T_c\subseteq\T'$.  The objects in
$\T'$ form a triangulated subcategory which is closed under taking
coproducts. It follows that $\T'=\T$ since $\T$ is compactly
generated. This finishes the proof.
\end{proof}

\section{Smashing subcategories}
Let $\S$ be a compactly generated triangulated category. In this
section we complete the correspondence between smashing localizations
of $\S$ and cohomological quotients of $\S_c$. In order to formulate
this, let us define the following full subcategories of $\S$ for any
ideal $\mathfrak I$ in $\S_c$:
\begin{align*}
\Filt\mathfrak I&= \{X\in\S\mid \mbox{every map $C\to X$, $C\in\S_c$,
factors through some map in $\mathfrak I$}\}, \\ \mathfrak I^\perp&=
\{X\in\S\mid \mbox{$\S(\p,X)=0$ for all $\p\in\mathfrak I$}\}.
\end{align*}

\begin{thm}\label{th:smash2}
Let $\S$ be a compactly generated triangulated category, and let
$\mathfrak I$ be an exact ideal in $\S_c$. Then there exists a
smashing localization $F\colon\S\to\T$ having the following
properties:
\begin{enumerate}
\item The right adjoint of $F$ identifies $\T$ with $\mathfrak I^\perp$.
\item $\Ker F=\Filt\mathfrak I$.
\item $\S_c\cap\Ann F=\mathfrak I$.
\end{enumerate}
\end{thm}

The proof of this result requires some preparations.  We start with
descriptions of $\Filt\mathfrak I$ and $\mathfrak I^\perp$ which we
take from \cite{K}.
\begin{lem}\label{le:ideal1}
Let $\mathfrak I$ be an ideal in $\S_c$ and $X$ be an object in $\S$.
\begin{enumerate}
\item $X\in\Filt\mathfrak I$ if and only if $H_X\in\li\Im\mathfrak I$. 
\item $X\in \mathfrak I^\perp$ if and only if
$\Hom_{\S_c}(\Im\mathfrak I,H_X)=0$.
\end{enumerate}
\end{lem}
\begin{proof}
For (1), see Lemma~3.9 in \cite{K}. (2) follows from the fact that
$\Hom_{\S_c}(-,H_X)$ is exact when restricted to $\mod\S_c$. 
\end{proof}

Now suppose that $\mathfrak I$ is a cohomological ideal in $\S_c$ and
observe that $\L=\li\Im\mathfrak I$ is a localizing subcategory of
$\Mod\S_c$, by Lemma~\ref{le:localizing}. Thus we obtain a quotient
functor $Q\colon\Mod\S_c\to\Mod\S_c/\L$ which has a right adjoint $R$;
see \cite[Proposition~III.8]{G}. Note that $R$ identifies
$\Mod\S_c/\L$ with the full subcategory $\M$ of $\S_c$-modules $M$
satisfying $\Hom_{\S_c}(\L,M)=0=\Ext^1_{\S_c}(\L,M)$; see
Lemma~\ref{le:abelian2}. Moreover, every $\S_c$-module $M$ fits into
an exact sequence
$$0\longrightarrow M'\longrightarrow M\longrightarrow (R\comp Q)M
\longrightarrow M''\longrightarrow 0$$
with $M',M''$ in $\L$.

\begin{lem}\label{le:ideal2} 
An object $X$ in $\S$ belongs to $\mathfrak I^\perp$ if and only if
$\Hom_{\S_c}(\L,H_X)=0=\Ext^1_{\S_c}(\L,H_X)$.
\end{lem}
\begin{proof} Suppose $X\in\mathfrak I^\perp$.
Then we have $\Hom_{\S_c}(\L,H_X)=0$ because $\Hom_{\S_c}(\Im\mathfrak
I,H_X)=0$ by Lemma~\ref{le:ideal1}.  Thus we have an exact sequence
$$0\longrightarrow H_X\longrightarrow (R\comp Q)H_X \longrightarrow
M\longrightarrow 0.$$ We claim that $M=0$. For this it is sufficient
to show that every map $\p\colon M'\to M$ from a finitely presented
module $M'$ is zero. The map $\p$ factors through some $M''$ in
$\Im\mathfrak I$ because $M\in\L$. Now we use that
$\Ext^1_{\S_c}(-,H_X)$ vanishes on finitely presented modules; see
\cite[Lemma~1.6]{K}.  Thus $M''\to M$ factors through $(R\comp
Q)H_X$. However, $\Hom_{\S_c}(\L,\M)=0$, and this implies $\p=0$.
Therefore $M=0$, and $H_X$ belongs to $\M$. Thus the proof is complete
because the other implication is trivial.
\end{proof}

\begin{lem}\label{le:perp}
Let $\mathfrak I$ be an ideal in $\S_c$ satisfying $\Si\mathfrak
I=\mathfrak I$. If $\mathfrak I$ is exact or $\mathfrak I^2=\mathfrak
I$, then $\mathfrak I^\perp$ is a triangulated subcategory of $\S$
which is closed under taking products and coproducts.
\end{lem}
\begin{proof}
Clearly, $\mathfrak I^\perp$ is closed under taking products and
coproducts. Also, $\Si(\mathfrak I^\perp)=\mathfrak I^\perp$ is
clear. It remains to show that $\mathfrak I^\perp$
is closed under forming extensions. Let
$$X\stackrel{\a}\lto Y\stackrel{\b}\lto
Z\stackrel{\g}\lto \Si X$$ be a triangle in $\S$ with $X$ and $Y$ in
$\mathfrak I^\perp$, which induces an exact sequence
$$0\lto\Coker H_\a\lto H_Z\lto\Ker H_{\Si\a}\lto 0$$
in $\Mod\S_c$.
The modules annihilated by $\mathfrak I$ form a subcategory which is
automatically closed under subobjects and quotients. The subcategory
is closed under extensions if $\mathfrak I^2=\mathfrak I$.  Thus $H_Z$
is annihilated by $\mathfrak I$ because $\Coker H_\a$ and $\Ker
H_{\Si\a}$ are annihilated by $\mathfrak I$. It follows that $Z$
belongs to $\mathfrak I^\perp$.  Now suppose that $\mathfrak I$ is
exact.  We apply Proposition~\ref{pr:exactideal} to see that the
category $\M$ of $\S_c$-modules $M$ satisfying
$\Hom_{\S_c}(\L,M)=0=\Ext^1_{\S_c}(\L,M)$ is closed under taking
kernels, cokernels, and extensions.  Thus $\Coker H_\a$ and $\Ker
H_{\Si\a}$ belong to $\M$, and therefore $H_Z$ as well.  We conclude
again that $Z$ belongs to $\mathfrak I^\perp$. This finishes the
proof.
\end{proof}

We are now in the position that we can prove Theorem~\ref{th:smash2}.

\begin{proof}[Proof of Theorem~\ref{th:smash2}.]
We know from Lemma~\ref{le:localizing} that $\L=\li\Im\mathfrak I$ is a localizing
subcategory of $\Mod\S_c$. Denote by $R$ the right adjoint of the
quotient functor $\Mod\S_c\to\Mod\S_c/\L$. Recall that $R$ identifies
$\Mod\S_c/\L$ with the full subcategory $\M$ of $\S_c$-modules $M$
satisfying $\Hom_{\S_c}(\L,M)=0=\Ext^1_{\S_c}(\L,M)$.  The quotient
$\Mod\S_c/\L$ is an abelian Grothendieck category. Thus there is an
injective cogenerator, say $I$, and we denote by $U$ the object in
$\S$ satisfying
$$\Hom_{\S_c}(H_\S-,RI)\cong\S(-,U).$$ Let $\T=\mathfrak I^\perp$,
which is a triangulated subcategory of $\S$ and closed under taking
products and coproducts by Lemma~\ref{le:perp}.  We claim that $U$ is
a perfect cogenerator for $\T$.  To see this, let $X$ be an object in
$\T$ satisfying $\S(X,U)=0$.  We have $H_X\in\M$ by
Lemma~\ref{le:ideal2}, and therefore $H_X=0$ since $RI$ is a
cogenerator for $\M$. Thus $X=0$. Now let $X_i\to Y_i$ be a set of
maps in $\T$ such that the map $\S(Y_i,U)\to\S(X_i,U)$ is surjective
for all $i$.  Using the fact that $RI$ is an injective cogenerator for
$\M$, we see that each map $H_{X_i}\to H_{Y_i}$ is a
monomorphism. Thus their product is a monomorphism, and therefore the
map $H_{\prod_iX_i}\to H_{\prod_iY_i}$ is a monomorphism. We conclude
that the map $\S(\prod_iY_i,U)\to\S(\prod_iX_i,U)$ is surjective,
since $RI$ is an injective object.  Now we can apply
Corollary~\ref{co:brown}.  It follows that the inclusion $\T\to\S$ has
a left adjoint $F$ which is a smashing localization since $\T$ is
closed under taking coproducts.

It remains to describe $\Ker F$ and $\Ann F$.  To this end let $X$ be
an object in $\S$.  We have $FX=0$ iff $\S(X,U)=0$ iff
$\Hom_{\S_c}(H_X,RI)=0$ iff $H_X\in\L$ iff $X\in\Filt\mathfrak I$, by
Lemma~\ref{le:ideal1}.  Now let $\p\colon X\to Y$ be a map in
$\S_c$. If $F\p=0$, then $\S(\p,U)=0$.  Now use that $\S(\p,U)=0$ iff
$\Hom_{\S_c}(H_\p,RI)=0$ iff $\Im H_\p\in\L$ iff $\Im
H_\p\in\Im\mathfrak I$ iff $\p\in\mathfrak I$.  Conversely, suppose
$\p\in\mathfrak I$. Then $FY\in\mathfrak I^\perp$ implies $F\p=0$
since $\T(F\p,FY)\cong\S(\p,FY)$. Thus $\S_c\cap\Ann F=\mathfrak I$,
and the proof is complete.
\end{proof}

Combining Theorem~\ref{th:smash1} and Theorem~\ref{th:smash2}, one
obtains a bijection between smashing localizations of $\S$ and exact
ideals of $\S_c$.  It is convenient to formulate this in terms of
smashing subcategories. Recall that a subcategory of $\S$ is {\em
smashing} if it of the form $\Ker F$ for some smashing localization
functor $F\colon\S\to\T$. Note that the kernel $\Ker F$ of any
localization functor $F$ is a {\em localizing} subcategory, that is,
$\Ker F$ is a full triangulated subcategory which is closed under
taking coproducts. Thus a subcategory $\R$ of $\S$ is smashing if and
only if $\R$ is a localizing subcategory admitting a right adjoint for
the inclusion $\R\to\S$ which preserves coproducts.

\begin{cor}\label{co:bijection}
Let $\S$ be a compactly generated triangulated category.
Then the maps 
$$\mathfrak I\mapsto \Filt\mathfrak I\quad\textrm{and}\quad
\R\mapsto\{\p\in\S_c\mid\mbox{$\p$ factors through some object in
$\R$}\}$$ induce mutually inverse bijections between the set of exact
ideals of $\S_c$ and the set of smashing subcategories of $\S$.
\end{cor}

A similar result appears as Theorem~4.9 in \cite{K}. However, the
proof given there is not correct for two reasons: it uses an
unnecessary assumption and relies on an erroneous definition of an
exact ideal.\begin{footnote}{The error occurs in Lemma~4.10. The claim
that $f^*(\mu)$ is an isomorphism is only correct if $f$ is a
cohomological quotient functor.}\end{footnote}

Let us formulate further consequences of Theorem~\ref{th:smash2} about
cohomological quotient functors.  I am grateful to B. Keller for
pointing out to me the following simple description of the exact
ideals of $\S_c$.

\begin{cor}\label{co:exact}
Let $\S$ be a compactly generated triangulated category.  Then an
ideal $\mathfrak I$ of $\S_c$ is exact if and only if the following
conditions hold.
\begin{enumerate}
\item $\Si\mathfrak I=\mathfrak I$.
\item $\mathfrak I$ is saturated.
\item $\mathfrak I$ is idempotent, that is, $\mathfrak I^2=\mathfrak I$.
\end{enumerate}
\end{cor}
\begin{proof}
Suppose first that $\mathfrak I$ is exact.  Applying
Corollary~\ref{co:bijection}, the ideal $\mathfrak I$ is the
collection of maps in $\S_c$ which factor through an object in
$\Filt\mathfrak I$.  Now fix a map $\p\colon X\to Y$ in $\mathfrak I$.
Then $\p$ factors through an object $Y'$ in $\Filt\mathfrak I$ via a
map $\p'\colon X\to Y'$, and $\p'$ factors through a map $\p_1\colon
X\to Y''$ in $\mathfrak I$ since $Y'$ belongs to $\Filt\mathfrak
I$. Thus $\p=\p_2\comp\p_1$, and $\p_2\colon Y''\to Y$ belongs to
$\mathfrak I$ because it factors through an object in $\Filt\mathfrak
I$.

Now suppose that $\mathfrak I$ satisfies (1) -- (3). The proof of
Theorem~\ref{th:smash2} works with these assumptions, thanks to
Lemma~\ref{le:perp}. The conclusion of Theorem~\ref{th:smash2} shows
that $\mathfrak I=\Ann F_c$ for some smashing localization
$F\colon\S\to\T$.  Thus $\mathfrak I$ is exact because $F_c$ is a
cohomological quotient functor by Theorem~\ref{th:smash1}.
\end{proof}

One may think of the following result as a generalization of the
localization theorem of Neeman-Ravenel-Thomason-Trobaugh-Yao
\cite{N2,R,TT,Y}.  To be precise, Neeman et al.\ considered
cohomological quotient functors of the form $\S_c\to\S_c/\R_0$ for
some triangulated subcategory $\R_0$ of $\S_c$ and analyzed the
smashing localization functor $\S\to\S/\R$ where $\R$ denotes the
localizing subcategory generated by $\R_0$.

\begin{cor}\label{co:adjoint}
Let $\S$ be a compactly generated triangulated category, and let 
$F\colon\S_c\to\T$ be a cohomological quotient functor. 
\begin{enumerate}
\item The category $\R=\Filt(\Ann F)$ is a smashing localizing subcategory of $\S$
and the quotient functor $\S\to\S/\R$ induces a
fully faithful and exact functor $\T\to\S/\R$ making the following
diagram commutative.
$$\xymatrix{\S_c\ar[rr]^F\ar[d]^{\mathrm{inc}}&&\T\ar[d]\\
\S\ar[rr]^{\mathrm{can}}&&\S/\R}$$
\item The triangulated category $\S/\R$ is compactly generated and
$(\S/\R)_c$ is the closure of the image of $\T\to\S/\R$ under forming
direct factors.
\item There exists a fully faithful and exact functor
$G\colon\T\to\S$ such that
$$\S(X,GY)\cong\T(FX,Y)$$
for all $X$ in $\S_c$ and $Y$ in $\T$.
\end{enumerate}
\end{cor}
\begin{proof}
The ideal $\Ann F$ is exact and we obtain from Theorem~\ref{th:smash2}
a smashing localization functor $Q\colon\S\to\S/\R$.  The induced
functor $Q_c\colon\S_c\to(\S/\R)_c$ is a cohomological quotient
functor with $\Ann Q_c=\Ann F$, by Theorem~\ref{th:smash1}. The proof
of Lemma~\ref{le:dominant} shows that $Q_c$ factors through $F$, since
idempotents in $(\S/\R)_c$ split. Moreover, the functor $\T\to\S/\R$
is fully faithful since it induces an equivalence
$\mod\T\to\mod(\S/\R)_c$.  Note that every compact object in $\S/\R$ is
a direct factor of some object in the image of $Q_c$ by
Theorem~\ref{th:cohquotient}. To obtain the functor $G\colon\T\to\S$,
take the fully faithful right adjoint $\S/\R\to \S$ of $Q$, and
compose this with the functor $\T\to\S/\R$.
\end{proof}

\section{The telescope conjecture}

The telescope conjecture due to Bousfield and Ravenel is originally
formulated for the stable homotopy category of CW-spectra; see
\cite[3.4]{B}, \cite[1.33]{R} (and \cite{MRS} for an unsuccessful
attempt to disprove the conjecture). The stable homotopy category is a
compactly generated triangulated category.  This fact suggests the
following formulation of the telescope conjecture for a specific
triangulated category $\S$ which is compactly generated.

\begin{telconj}
Every smashing subcategory of $\S$ is generated as a localizing
subcategory by objects which are compact in $\S$.
\end{telconj}

Recall that a subcategory of $\S$ is {\em smashing} if it is of the form
$\Ker F$ for some smashing localization functor $F\colon\S\to\T$. Note
that $\Ker F$ is a {\em localizing} subcategory of $\S$, that is,
$\Ker F$ is a full triangulated subcategory which is closed under
taking coproducts. A localizing subcategory of $\S$ is {\em generated}
by a class $\X$ of objects if it is the smallest localizing
subcategory of $\S$ which contains $\X$.

The telescope conjecture in this general form is known to be false. In
fact, Keller gives an example of a smashing subcategory which contains
no non-zero compact object \cite{Ke}; see also
Section~\ref{se:hloc}. However, there are classes of compactly
generated triangulated categories where the conjecture has been
verified.

We have seen that smashing subcategories of $\S$ are closely related
to cohomological quotients of $\S_c$. It is therefore natural to
translate the telescope conjecture into a statement about
cohomological quotients. Roughly speaking, the telescope conjecture
for $\S$ is equivalent to the assertion that every flat epimorphism
$\S_c\to\T$ is an Ore localization. We need some preparations in order
to make this precise.

\begin{lem}\label{le:tele}
Let $\S$ be a compactly generated triangulated category and
$F\colon\S\to\T$ be a smashing localization functor.  Then the
following are equivalent.
\begin{enumerate}
\item The localizing subcategory $\Ker F$ is generated by objects
which are compact in $\S$.
\item The ideal $\S_c\cap \Ann F$ of $\S_c$ is generated by identity maps. 
\end{enumerate}
\end{lem}
\begin{proof}
Let $\mathcal R=\Ker F$ and $\mathfrak I=\S_c\cap \Ann F$.

(1) $\Rightarrow$ (2): Suppose $\R$ is generated by compact
objects. Then every object in $\R$ is a homotopy colimit of objects in
$\R\cap\S_c$; see \cite{N1}. Let $\p\colon X\to Y$ be a map in
$\mathfrak I$. It follows from Lemma~\ref{le:ann} that $\p$ factors
through a homotopy colimit of objects in $\R\cap\S_c$. Thus $\p$
factors through some object in $\R\cap\S_c$ since $X$ is compact. We
conclude that $\mathfrak I$ is generated by the identity maps of all
objects in $\R\cap\S_c$.

(2) $\Rightarrow$ (1): Let $\R_0$ be a class of compact objects and
suppose $\mathfrak I$ is generated by the identity maps of all objects
in $\R_0$.  Let $\R'$ be the localizing subcategory which is generated
by $\R_0$. This category is smashing and we have a localization
functor $F'\colon\S\to\T'$ with $\R'=\Ker F'$. Clearly, $\Ann
F_c\subseteq \Ann F'_c$ since $\id_X$ belongs to $\Ann F'_c$ for all
$X$ in $\R_0$. On the other hand, $\R'\subseteq \R$ since
$\R_0\subseteq\R$.  Thus $\Ann F'_c\subseteq \Ann F_c$ by
Lemma~\ref{le:ann}, and $\Ann F'_c= \Ann F_c$ follows.  We conclude
that $\R'=\R$ because Theorem~\ref{th:smash2} states that a smashing
subcategory is determined by the corresponding exact ideal in $\S_c$.
Thus $\R$ is generated by compact objects.
\end{proof}

\begin{prop}\label{pr:tele}
Let $F\colon\C\to\D$ be an exact  functor between triangulated categories.
Then the following are equivalent.
\begin{enumerate}
\item $F$ induces an equivalence $\C/{\Ker F}\to\D$.
\item $F$ induces an equivalence $\C[\Phi^{-1}]\to\D$ where
$\Phi=\{\p\in\C\mid \mbox{$F\p$ is an iso}\}$.
\item $F$ is a CE-quotient functor and the ideal $\Ann F$ is generated
by identity maps.
\end{enumerate}
\end{prop}
\begin{proof} 
We put $\B=\Ker F$ and denote by $Q\colon\C\to\C/\B$ the quotient
functor.

(1) $\Leftrightarrow$ (2): The quotient $\C/\B$ is by definition
$\C[\Psi^{-1}]$ where $\Psi$ is the class of maps $X\to Y$ in $\C$
which fit into an exact triangle $X\to Y\to Z\to\Si X$ with $Z$ in
$\B$. The exactness of $F$ implies that $\Psi$ is precisely the class
of maps $\p$ in $\C$ such that $F\p$ is invertible.

(1) $\Rightarrow$ (3): We have seen in Example~\ref{ex:verdier1} and
Example~\ref{ex:verdier2} that $Q$ is a CE-quotient
functor. Lemma~\ref{le:ann} implies that the ideal $\Ann Q$
is generated by the identity maps of all objects in $\B$.

(3) $\Rightarrow$ (1): The functor $F$ induces an exact functor
$\C/\B\to\D$. Now suppose that $\Ann F$ is generated by identity
maps. Then $\Ann F$ is generated by the identity maps of all objects
in $\B$, and therefore $\Ann Q=\Ann F$ by Lemma~\ref{le:ann}.  Thus
$Q$ factors through $F$ by an exact functor $\D\to\C/\B$ since $F$ is
an exact quotient functor. The uniqueness of $\D\to\C/\B$ and
$\C/\B\to\D$ implies that both functors are mutually inverse
equivalences.
\end{proof}

\begin{lem}\label{le:annideal} 
Let $F\colon\C\to\D$ be a cohomological quotient functor.  Then the
following are equivalent.
\begin{enumerate}
\item $F$ induces a fully faithful functor $\C/{\Ker F}\to\D$.
\item The ideal $\Ann F$ is generated by identity maps.
\end{enumerate}
\end{lem}
\begin{proof}
Let $\D'$ be the smallest full triangulated subcategory of $\D$
containing the image of $F$. It follows from Lemma~\ref{le:dominant}
that the induced functor $F'\colon\C\to\D'$ is a CE-quotient functor.
Now the assertion follows from Proposition~\ref{pr:tele} since $\Ann
F=\Ann F'$.
\end{proof}

We obtain the following reformulation of the telescope
conjecture. Note in particular, that the telescope conjecture becomes
a statement about the category of compact objects.

\begin{thm}\label{th:tele}
Let $\S$ be a compactly generated triangulated category. Then the
following are equivalent.
\begin{enumerate}
\item Every smashing subcategory of $\S$ is generated by objects which
are compact in $\S$.
\item Every smashing subcategory of $\S$ is a compactly generated
triangulated category.
\item Every exact ideal in $\S_c$ is generated by idempotent elements.
\item Every CE-quotient functor $F\colon \S_c\to\T$ induces an
equivalence $\S_c/{\Ker F}\to\T$.
\item Every cohomological quotient functor $F\colon \S_c\to\T$ induces
a fully faithful functor $\S_c/{\Ker F}\to\T$.
\item Every flat epimorphism $F\colon \S_c\to\T$ satisfying $\Si(\Ann
F)=\Ann F$ induces an equivalence $\S_c[\Phi^{-1}]\to\T$ where
$\Phi=\{\p\in\S_c\mid \mbox{$F\p$ is an iso}\}$.
\end{enumerate}
\end{thm}
\begin{proof} 
We use the bijection between smashing subcategories of $\S$ and exact
ideals of $\S_c$; see Corollary~\ref{co:bijection}. Recall that an
ideal is by definition exact if it is of the form $\Ann F$ for some
cohomological functor $F\colon\S_c\to\T$.

(1) $\Leftrightarrow$ (2): The inclusion $\R\to\S$ of a smashing
subcategory preserves compactness.

(1) $\Leftrightarrow$ (3): Apply Lemma~\ref{le:tele}.  Note that any
ideal in $\S_c$ which is generated by idempotent maps is also
generated by identity maps. This follows from the fact that
idempotents in $\S_c$ split.

(3) $\Leftrightarrow$ (4): Apply Proposition~\ref{pr:tele}.

(3) $\Leftrightarrow$ (5): Apply Lemma~\ref{le:annideal}.

(5) $\Rightarrow$ (6): Let $F\colon \S_c\to\T$ be a flat epimorphism
satisfying $\Si(\Ann F)=\Ann F$. Composing it with the idempotent
completion $\T\to\bar\T$ gives a cohomological quotient functor
$\S_c\to\bar\T$, by Theorem~\ref{th:cohquotient} and
Theorem~\ref{th:epi}. Now use that $\S_c[\Phi^{-1}]=\S_c/{\Ker
F}$. Thus $\S_c[\Phi^{-1}]\to\T$ is fully faithful, and it is an
equivalence since $F$ is surjective on objects by Lemma~\ref{le:epi2}.

(6) $\Rightarrow$ (5): Let $F\colon \S_c\to\T$ be a cohomological
quotient functor, and denote by $\T'$ the full subcategory of $\T$
whose objects are those in the image of $F$. The induced functor
$\S_c\to\T'$ is a flat epimorphism. Now use again that
$\S_c[\Phi^{-1}]=\S_c/{\Ker F}$. Thus the induced functor $\S_c/{\Ker
F}\to\T$ is fully faithful.
\end{proof}

\begin{rem} 
Let $\C$ be a triangulated category and $F\colon\C\to\D$ be a flat
functor satisfying $\Si(\Ann F)=\Ann F$. Then $\Phi=\{\p\in\C\mid
\mbox{$F\p$ is an iso}\}$ is a {\em multiplicative system}, that is,
$\Phi$ admits a calculus of left and right fractions in the sense of
\cite{GZ}.
\end{rem}

We say that an additive functor $\C\to\D$ is an {\em Ore localization}
if it induces an equivalence $\C[\Phi^{-1}]\to\D$ for some
multiplicative system $\Phi$ in $\C$.  Using this terminology,
Theorem~\ref{th:tele} suggests the following reformulation of the
telescope conjecture.

\begin{cor} 
The telescope conjecture holds true for a compactly generated
triangulated category $\S$ if and only if every flat epimorphism
$F\colon \S_c\to\T$ satisfying $\Si(\Ann F)=\Ann F$ is an Ore
localization.
\end{cor}

The reformulation of the telescope conjecture in terms of exact ideals
raises the question when an idempotent ideal is generated by
idempotent elements.  This follows from Corollary~\ref{co:exact} where
it is shown that the exact ideals are precisely the idempotent ideals
which satisfy some natural extra conditions. 

\begin{cor} 
The telescope conjecture holds true for a compactly generated
triangulated category $\S$ if and only if every idempotent and
saturated ideal $\mathfrak I$ of $\S_c$ satisfying $\Si\mathfrak
I=\mathfrak I$ is generated by idempotent elements.
\end{cor}

The problem of finding idempotent generators for an idempotent ideal
is a very classical one from ring theory. For instance, Kaplansky
introduced the class of SBI rings, where SBI stands for `suitable for
building idempotent elements' \cite[III.8]{J}. Also, Auslander asked
the question for which rings every idempotent ideal is generated by an
idempotent element \cite[p.\ 241]{A}.  One can show for an additive
category $\C$, that every idempotent ideal is generated by idempotent
elements provided that $\C$ is perfect in the sense of Bass
\cite{Mi}. Recall that $\C$ is {\em perfect} if every object in $\C$
decomposes into a finite coproduct of indecomposable objects with
local endomorphism rings, and, for every sequence
$$X_1\stackrel{\p_1}\lto X_2\stackrel{\p_2}\lto
X_3\stackrel{\p_3}\lto\cdots$$ of non-isomorphisms between
indecomposable objects, the composition
$\p_n\comp\ldots\comp\p_2\comp\p_1$ is zero for $n$ sufficiently
large. Note that the category $\S_c$ of compact objects is perfect if
and only if every object in $\S$ is a coproduct of indecomposable
objects with local endomorphism rings \cite[Theorem~2.10]{K}.

\section{Homological epimorphisms of rings}\label{se:hepi}

A commutative localization $R\to S$ of an associative ring $R$ is
always a flat epimorphism. For non-commutative localizations, there is
a weaker condition which is often satisfied. Recall from \cite{GL}
that a ring homomorphism $R\to S$ is a {\em homological epimorphism}
if $S\otimes_RS\cong S$ and $\Tor_i^R(S,S)=0$ for all $i\geq
1$. Homological epimorphisms frequently arise in representation theory
of finite dimensional algebras, in particular via universal
localizations \cite{C,S}.

In recent work of Neeman and Ranicki \cite{NR}, homological
epimorphisms appear when they study the following chain complex lifting
problem for a ring homomorphism $R\to S$.  We denote by $\mathbf
K^b(R)$ the homotopy category of bounded complexes of finitely
generated projective $R$-modules.

\begin{defn}
Fix a ring homomorphism $R\to S$.
\begin{enumerate}
\item We say that the {\em chain complex lifting problem} has a
positive solution, if every complex $Y$ in $\mathbf K^b(S)$ such that
for each $i$ we have $Y^i=P^i\otimes_RS$ for some finitely generated
projective $R$-module $P^i$, is isomorphic to $X\otimes_RS$ for some
complex $X$ in $\mathbf K^b(R)$.
\item We say that the {\em chain map lifting problem} has a positive
solution, if for every pair $X,Y$ of complexes in $\mathbf K^b(R)$ and
every map $\a\colon X\otimes_RS\to Y\otimes_RS$ in $\mathbf K^b(S)$,
there are maps $\phi\colon X'\to X$ and $\a'\colon X'\to Y$ in
$\mathbf K^b(R)$ such that $\phi\otimes_RS$ is invertible and
$\a=\a'\otimes_RS\comp(\phi\otimes_RS)^{-1}$ in $\mathbf K^b(S)$.
\end{enumerate}
\end{defn}

The following observation shows that both lifting problems are closely
related. In fact, it seems that the more general problem of lifting
maps is the more natural one. 

\begin{lem}\label{le:cclp}
Given a ring homomorphism $R\to S$, the chain complex lifting problem
has a positive solution whenever the chain map lifting problem has a
positive solution.
\end{lem}
\begin{proof}
Fix a complex $Y$ in $\mathbf K^b(S)$. We proceed by induction on its
length $\ell(Y)=n$. If $n=0$, then $Y$ is concentrated in one degree,
say $i$, and therefore $Y=X\otimes_RS$ for $X=P^i$. If $n>0$, choose
an exact triangle $Y_1\to Y_2\to Y\to\Si Y_1$ with $\ell(Y_i)<n$ for
$i=1,2$.  By our assumption, we have $Y_i\cong X_i\otimes_RS$ for some
complexes $X_1$ and $X_2$ in $\mathbf K^b(R)$. Moreover, using the
positive solution of the chain map lifting problem, the map
$X_1\otimes_RS\to X_2\otimes_RS$ is of the form
$\a\otimes_RS\comp(\phi\otimes_RS)^{-1}$ for some maps $\phi\colon
X_1'\to X_1$ and $\a\colon X_1'\to X_2$ in $\mathbf K^b(R)$.  We
complete $\a$ to an exact triangle $X_1'\to X_2\to X\to\Si X_1'$ and
conclude that $Y\cong X\otimes_RS$.
\end{proof}

The example of a proper field extension $k\to K$ shows that both
lifting problems are not equivalent. In fact, the chain complex
lifting problem for $k\to K$ has a positive solution, but the chain
map lifting problem does not.

The proof of Lemma~\ref{le:cclp} suggests the following reformulation
of the chain complex lifting problem.

\begin{lem}
Given a ring homomorphism $R\to S$, the chain complex lifting problem
has a positive solution if and only if the full subcategory of
$\mathbf K^b(S)$ formed by the objects in the image of $-\otimes_RS$
is a triangulated subcategory.
\end{lem}
\begin{proof} Clear.
\end{proof}

We continue with a reformulation of the chain map lifting
problem.

\begin{prop}\label{pr:maplifting}
Let $R\to S$ be a ring homomorphism and denote by $\K$ the full
subcategory of $\mathbf K^b(R)$ formed by the complexes $X$ such that
$X\otimes_RS=0$.  Then the following are equivalent.
\begin{enumerate}
\item The chain map lifting problem has a positive solution.
\item The functor $-\otimes_RS$ induces a fully faithful functor $\mathbf
K^b(R)/\K\to\mathbf K^b(S)$.
\end{enumerate}
\end{prop}
\begin{proof}
(1) $\Rightarrow$ (2): Denote by $\D$ the full subcategory of $\mathbf
K^b(S)$ which is formed by all objects in the image of $-\otimes_RS$.
Using the description of the maps in $\D$, we observe that $\D$ is a
triangulated subcategory of $\mathbf K^b(S)$. It follows from
Proposition~\ref{pr:verdier} that $F$ induces an equivalence $\mathbf
K^b(R)/\K\to\D$.

(2) $\Rightarrow$ (1): The maps in $\mathbf K^b(R)/\K$ can be
described as fractions; see for instance Proposition~\ref{pr:verdier}.
The functor $\mathbf K^b(R)/\K\to\mathbf K^b(S)$ being full
implies the positive solution of the chain map lifting problem.
\end{proof}

Given a ring homomorphism $R\to S$, we shall see that the problem of
lifting complexes and their maps is closely related to the question,
when the derived functor $-\otimes_R^{\mathbf L}S$ is a smashing
localization.

Let us denote by $\mathbf D(R)$ the unbounded derived category of
$R$. Note that the inclusion $\mathbf K^b(R)\to\mathbf D(R)$ induces
an equivalence $\mathbf K^b(R)\to\mathbf D(R)_c$.

\begin{thm}\label{th:hepi}
For a ring homomorphism $R\to S$ the following are equivalent.
\begin{enumerate}
\item The derived functor $-\otimes_R^{\mathbf L}S\colon\mathbf
D(R)\to \mathbf D(S)$ is a smashing localization.
\item The functor $-\otimes_RS\colon\mathbf K^b(R)\to \mathbf K^b(S)$
is a cohomological quotient functor.
\item The map $R\to S$ is a homological epimorphism.
\end{enumerate}
\end{thm}
\begin{proof} The functor  $F=-\otimes_R^{\mathbf L}S\colon\mathbf
D(R)\to \mathbf D(S)$ has a right adjoint $G\colon\mathbf D(S)\to
\mathbf D(R)$ which is simply restriction of scalars, that is,
$G=\mathbf R\Hom_S(S,-)$.  Clearly, $G$ preserves coproducts. Thus $F$
is a smashing localization if and only if $F$ is a localization
functor. Note that $F$ is a localization functor if and only if
$F\comp G\cong\Id_{\mathbf D(S)}$.  Moreover, $F\comp G$ is exact and
preserves coproducts. Using infinite devissage, one sees that $F\comp
G\cong\Id_{\mathbf D(S)}$ if and only if the canonical map
$X\otimes_R^{\mathbf L}S\to X$ is an isomorphism for the complex $X=S$
which is concentrated in degree $0$. Clearly, this condition is
equivalent to $S\otimes_RS\cong S$ and $\Tor_i^R(S,S)=0$ for all
$i\geq 1$. This proves the equivalence of (1) and (3).  The
equivalence of (1) and (2) follows from Theorem~\ref{th:smash1} since
$F|_{\mathbf K^b(R)}=-\otimes_RS$.
\end{proof}

We obtain the following conditions for solving the chain map lifting
problem.

\begin{thm}\label{th:maplifting}
Given a ring homomorphism $R\to S$, the chain map lifting problem has
a positive solution if and only if
\begin{enumerate}
\item $R\to S$ is a homological epimorphism, and
\item every map $\p$ in $\mathbf K^b(R)$ satisfying $\p\otimes_RS=0$
factors through some $X$ in $\mathbf K^b(R)$ such that
$X\otimes_RS=0$.
\end{enumerate}
\end{thm}
\begin{proof} 
Suppose first that (1) and (2) hold.  Condition (1) says that
$F=-\otimes_RS\colon \mathbf K^b(R)\to\mathbf K^b(S)$ is a
cohomological quotient functor. This follows from
Theorem~\ref{th:hepi}. Applying Lemma~\ref{le:annideal}, we conclude
from (2) that $F$ induces a fully faithful functor $\mathbf
K^b(R)/{\Ker F}\to\mathbf K^b(S)$.  The positive solution of the chain
map lifting problem follows from Proposition~\ref{pr:maplifting}.

Now suppose that the chain map lifting problem has a positive
solution. The description of the maps in the image $\Im F$ of $F$
implies that the full subcategory formed by the objects in $\Im F$ is
a triangulated subcategory of $\mathbf K^b(S)$. It contains a
generator of $\mathbf K^b(S)$ and therefore every object in $\mathbf
K^b(S)$ is a direct factor of some object in $\Im F$.  Now we apply
Theorem~\ref{th:cohquotient} and see that $F$ is a cohomological
quotient functor. Thus (1) holds by Theorem~\ref{th:hepi}.  The
induced functor $\mathbf K^b(R)/{\Ker F}\to\mathbf K^b(S)$ is fully
faithful by Proposition~\ref{pr:maplifting}.  It follows from
Lemma~\ref{le:annideal} that (2) holds. This finishes the proof.
\end{proof}

\begin{cor}
Let $R$ be a ring such that the telescope conjecture holds true for
$\mathbf D(R)$. Then the chain map lifting problem has a positive
solution for a ring homomorphism $f\colon R\to S$ if and only if $f$
is a homological epimorphism.
\end{cor}

Note that the telescope conjecture has been verified for $\mathbf
D(R)$ provided $R$ is commutative noetherian \cite{N}. On the other
hand, Keller has given an example of a ring $R$ such that the
telescope conjecture for $\mathbf D(R)$ does not hold \cite{Ke}. Let
us mention that the validity of the telescope conjecture is preserved
under homological epimorphims.

\begin{prop} 
  Let $R\to S$ be a homological epimorphism. If the telescope
  conjecture holds for $\mathbf D(R)$, then the telescope conjecture
  holds for $\mathbf D(S)$.
\end{prop}
\begin{proof}
The derived functor $F=-\otimes_R^{\mathbf L}S\colon\mathbf D(R)\to
\mathbf D(S)$ is a smashing localization by Theorem~\ref{th:hepi}.
Now suppose that $G\colon\mathbf D(S)\to \T$ is a smashing
localization.  A composite of smashing localizations is a smashing
localizations.  Thus $\Ker F$ is generated by a class $\X$ of compact
objects since the telescope conjecture holds for $\mathbf D(R)$.
It follows that $\Ker G$ is generated by $F\X$.
\end{proof}

The work of Neeman and Ranicki \cite{NR} on the problem of lifting
chain complexes is motivated by some applications in algebraic
$K$-theory. In fact, they generalize the classical long exact sequence
which is induced by an injective Ore localization. More precisely,
they show that every universal localization $f\colon R\to S$ which is
a homological epimorphism induces a long exact sequence
$$\cdots \lto K_1(R)\lto K_1(S)\lto K_0(R,f)\lto K_0(R)\lto K_0(S)$$
in algebraic $K$-theory \cite[Theorem~10.11]{NR}. Our analysis of the
chain map lifting problem suggests that being a homological
epimorphism and satisfying the additional hypothesis (2) in
Theorem~\ref{th:maplifting} is the crucial property for such a
sequence.  We sketch the construction of this sequence which uses the
machinery developed by Waldhausen in \cite{Wa}. Our exposition follows
closely the ideas of Thomason-Trobaugh \cite{TT} and Neeman-Ranicki
\cite{NR}.

We fix a ring homomorphism $f\colon R\to S$. Denote by $\mathbf W(R)$
the complicial biWaldhausen category of bounded chain complexes of
finitely generated projective $R$-modules \cite[1.2.11]{TT}. We denote
by $K(R)$ the corresponding Waldhausen $K$-theory spectrum $K(\mathbf
W(R))$; see \cite[1.5.2]{TT}. Note that $K(R)$ is homotopy equivalent
to the Quillen $K$-theory spectrum of the exact category $\proj R$ of
finitely generated projective $R$-modules \cite[1.11.2]{TT}.  The
algebraic $K$-groups $K_n(R)=\pi_n K(R)$ are by definition the
homotopy groups of the spectrum $K(R)$. Now let $\mathbf W(R,f)$ be
the complicial biWaldhausen subcategory of $\mathbf W(R)$ consisting
of those complexes $X$ in $\mathbf W(R)$ such that $X\otimes_RS$ is
acyclic, and put $K(R,f)= K(\mathbf W(R,f))$.

\begin{thm}
Let $f\colon R\to S$ be a homological epimorphism and suppose $f$
satisfies condition (2) in Theorem~\ref{th:maplifting}.  Then $f$
induces a sequence
$$\mathbf W(R,f)\lto\mathbf W(R)\lto\mathbf W(S)$$
of exact functors such that
$$K(R,f)\lto K(R)\lto K(S)$$
is a homotopy fibre sequence, up to
failure of surjectivity of $K_0(R)\to K_0(S)$.  In particular, there
is induced a long exact sequence
$$\cdots \lto K_1(R)\lto K_1(S)\lto K_0(R,f)\lto K_0(R)\lto K_0(S)$$
of algebraic $K$-groups.
\end{thm}
\begin{proof}
  The proof is modeled after that of Thomason-Trobaugh's localization
  theorem \cite[Theorem~5.1]{TT}. We recall that a complicial
  biWaldhausen category comes equipped with cofibrations and weak
  equivalences \cite[1.2.11]{TT}.  The cofibrations of $\mathbf W(R)$
  are by definition the chain maps which are split monomorphism in
  each degree, and the weak equivalences are the quasi-isomorphisms.
  We define a new complicial biWaldhausen category $\mathbf W(R/f)$ as
  follows. The underlying category is that of $\mathbf W(R)$, the
  cofibrations are those of $\mathbf W(R)$, and the weak equivalences
  are the chain maps whose mapping cone lies in $\mathbf W(R,f)$. We
  denote by $K(R/f)$ the $K$-theory spectrum of $\mathbf W(R/f)$ and
  obtain an induced sequence
$$\mathbf W(R,f)\lto\mathbf W(R)\lto\mathbf W(R/f)$$
of exact functors such that
$$K(R,f)\lto K(R)\lto K(R/f)$$ is a homotopy fibre sequence by
Waldhausen's localization theorem \cite[1.8.2]{TT}.  The functor
$\mathbf W(R)\to\mathbf W(S)$ factors through $\mathbf W(R)\to\mathbf
W(R/f)$ and induces an exact functor $\mathbf W(R/f)\to\mathbf
W(S)$. Note that any exact functor $\mathbf A\to\mathbf B$ between
complicial biWaldhausen categories induces a homotopy equivalence of
$K$-theory spectra $K(\mathbf A)\to K(\mathbf B)$ provided the functor
induces an equivalence $\Ho(\mathbf A)\to\Ho(\mathbf B)$ of the
derived homotopy categories \cite[1.9.8]{TT}. Observe that
$\Ho(\mathbf W(R))=\mathbf K^b(R)$. Moreover, $\mathbf W(R)\to\mathbf
W(R/f)$ induces an equivalence
$$\Ho(\mathbf W(R))/{\Ho(\mathbf
W(R,f))}\stackrel{\sim}\lto\Ho(\mathbf W(R/f)).$$ Thus we have the
following commutative diagram
$$\xymatrix{ \Ho(\mathbf W(R,f))\ar[d]^\wr\ar[r]&\Ho(\mathbf
W(R))\ar[d]^\wr\ar[r] &\Ho(\mathbf W(R/f))\ar[d]^\wr\ar[r]&\Ho(\mathbf
W(S))\ar[d]^\wr\\ \K\ar[r]&\mathbf K^b(R)\ar[r]&\mathbf
K^b(R)/\K\ar[r] &\mathbf K^b(S)}$$ where $\K$ denotes the full
subcategory of $\mathbf K^b(R)$ consisting of all complexes $X$ such
that $X\otimes_RS=0$.  Next we use our assumption about the ring
homomorphism $f$ and apply Proposition~\ref{pr:maplifting} and
Theorem~\ref{th:maplifting}.  It follows that $\mathbf
W(R/f)\to\mathbf W(S)$ induces a functor
$$\Ho(\mathbf W(R/f))\lto\Ho(\mathbf W(S))$$ which is an equivalence
up to direct factors. We conclude from the cofinality theorem
\cite[1.10.1]{TT} that
$$K(R,f)\lto K(R)\lto K(S)$$ is a homotopy fibre sequence,
up to failure of surjectivity of $K_0(R)\to K_0(S)$.
\end{proof}

\section{Homological localizations of rings}\label{se:hloc}

Let $R$ be an associative ring and let $\Phi$ be a class of maps
between finitely generated projective $R$-modules. The {\em universal
localization} of $R$ with respect to $\Phi$ is the universal ring
homomorphism $R\to S$ such that $\p\otimes_RS$ is an isomorphism of
$S$-modules for all $\p$ in $\Phi$; see \cite{C,S}. To
construct $S$, one formally inverts all maps from $\Phi$ in the
category $\C=\proj R$ of finitely generated projective $R$-modules and
puts $S=\C[\Phi^{-1}](R,R)$. The following concept generalizes
universal localizations.

\begin{defn}
We call a ring homomorphism $f\colon R\to S$ a {\em homological
localization} with respect to a class $\Phi$ of complexes in $\mathbf
K^b(R)$ if
\begin{enumerate}
\item $X\otimes_RS=0$ in $\mathbf K^b(S)$ for all $X$ in $\Phi$, and
\item given any ring homomorphism $f'\colon R\to S'$ such that
$X\otimes_RS'=0$ in $\mathbf
K^b(S')$ for all $X$ in $\Phi$, there exists a unique
homomorphism $g\colon R\to R'$ such that $f'=g\comp f$.
\end{enumerate}
\end{defn}

Any universal localization is a homological localization.  In fact,
any map $\p\colon P\to Q$ between finitely generated projective
$R$-modules may be viewed as a complex of length one by taking its
mapping cone $\Cone\p$.  If $R\to S$ is a ring homomorphism, then
$\p\otimes_RS$ is an isomorphism if and only if $(\Cone\p)\otimes_RS=0$ in
$\mathbf K^b(S)$.

The following example, which I learned from A. Neeman, shows that a
homological localization need not to exist.

\begin{exm} 
  Let $k$ be a field and $R=k[x,y]$. Let $P$ be the complex $$
  \cdots
  \lto 0\lto R
  \stackrel{\left[\begin{smallmatrix}x\\y\end{smallmatrix}\right]}\lto
  R\amalg R
  \stackrel{\left[\begin{smallmatrix}y&-x\end{smallmatrix}\right]}\lto
  R \lto 0\lto\cdots$$
  which is a projective resolution of $R/(x,y)$.
  Then we have $P\otimes_RR[x^{-1}]=0$ and $P\otimes_RR[y^{-1}]=0$.
  Now suppose there is a homological localization $R\to S$ with
  respect to $P$.  Viewing $R[x^{-1}]$ and $R[y^{-1}]$ as subrings of
  $k(x,y)$, we have $R[x^{-1}]\cap R[y^{-1}]=R$.  Therefore the
  identity $R\to R$ factors through $R\to S$. This is a contradiction
  and shows that the homological localization with respect to $P$
  cannot exist.
\end{exm}

Next we consider an example of a homological epimorphism which is not
a homological localization. Keller used this example in order to
disprove the telescope conjecture for the derived category of a ring
\cite{Ke}.  Let us explain the idea of Keller's example. He uses the
following observation.

\begin{lem}\label{le:wod}
Let $R$ be a ring and $\mathfrak a$ be a two-sided ideal which is
contained in the Jacobson radical of $R$. Then $X\otimes_RR/\mathfrak
a=0$ implies $X=0$ for every $X$ in $\mathbf K^b(R)$.
\end{lem}
\begin{proof} 
Using induction on the length of the complex $X$, the assertion
follows from Nakayama's lemma.
\end{proof}

In \cite{W}, Wodzicki has constructed an example of a ring $R$ such
that the Jacobson radical $\mathfrak r$ is non-zero and satisfies
$$\Tor_i^R(R/{\mathfrak r},R/{\mathfrak r})=0\quad\textrm{for
all}\quad i\geq 1.$$ Thus $R\to R/{\mathfrak r}$ is a homological
epimorphism which induces a cohomological quotient functor
$$F=-\otimes_RR/{\mathfrak r}\colon \mathbf K^b(R)\lto\mathbf
K^b(R/{\mathfrak r})$$
satisfying $\Ker F=0$ and $\Ann F\neq 0$. It
follows that $R\to R/{\mathfrak r}$ is not a homological localization
since $\Ker F=0$. Moreover, Theorem~\ref{th:tele} shows that the
telescope conjecture does not hold for $\mathbf D(R)$.

One can find more examples along these lines, as R. Buchweitz kindly
pointed out to me. Take any {\em B\'ezout domain} $R$, that is, an
integral domain such that every finitely generated ideal is principal.
We have for every ideal $\mathfrak a$
$$\Tor_1^R(R/\mathfrak a,R/\mathfrak a)\cong\mathfrak a/\mathfrak a^2
\quad\textrm{and}\quad\Tor_i^R(R/\mathfrak a,-)=0 \quad\textrm{for
all}\quad i> 1.$$ Thus for any idempotent ideal $\mathfrak a$, the
natural map $R\to R/\mathfrak a$ is a homological epimorphism.
Specific examples arise from valuation domains, which are precisely
the local B\'ezout domains.

Our interest in homological localizations is motivated by the
following observation, which shows that the positive solution of the
chain map lifting problem forces a ring homomorphism to be a
homological localization.

\begin{prop} 
Let $f\colon R\to S$ be a ring homomorphism and suppose the chain map
lifting problem has a positive solution. Then $f$ is a homological
localization.
\end{prop}
\begin{proof}
Denote by $\Phi$ the set of complexes $X$ in $\mathbf K^b(R)$ such
that $X\otimes_RS=0$, and denote by $\K$ the corresponding full
subcategory.  We have seen in Proposition~\ref{pr:maplifting} that
$-\otimes_RS\colon\mathbf K^b(R)\to\mathbf K^b(S)$ induces a fully
faithful functor $\mathbf K^b(R)/\K\to\mathbf K^b(S)$.  Now suppose
that $f'\colon R\to S'$ is a ring homomorphism satisfying
$X\otimes_RS'=0$ for all $X$ in $\Phi$. Then $-\otimes_RS'$ factors
through the quotient functor $\mathbf K^b(R)\to\mathbf K^b(R)/\K$ via
some functor $G\colon \mathbf K^b(R)/\K\to\mathbf K^b(S')$.  Clearly,
$G$ induces a homomorphism $g\colon S\to S'$ such that $f'=g\comp f$.
The uniqueness of $g$ follows from the uniqueness of $G$.
\end{proof}

In \cite{NR}, Neeman and Ranicki show that the chain complex lifting
problem has a positive solution for every universal localization which
is a homological epimorphism. We give an alternative proof of this
result which is based on the criterion for lifting chain maps in
Theorem~\ref{th:maplifting}.

\begin{thm}\label{th:universal}
The chain map lifting problem has a positive solution for every
homological epimorphism $R\to S$ which is a  universal localization.
\end{thm}

We need some preparations for the proof of this result. Fix a
homological epimorphism $R\to S$, and suppose it is the universal
localization with respect to a class $\Phi$ of maps in the category
$\C=\proj R$ of finitely generated projective $R$-modules.  Thus we
have $\proj S=\C[\Phi^{-1}]$. We denote by
$\Cone\Phi=\{\Cone\p\mid\p\in\Phi\}$ the corresponding class of complexes
of length one in $\mathbf K^b(R)=\mathbf K^b(\C)$, and we write
$\langle\Cone\Phi\rangle$ for the thick subcategory generated by
$\Cone\Phi$. Finally, denote by $\T$ the idempotent completion of the
quotient $\mathbf K^b(\C)/{\langle\Cone\Phi\rangle}$, which one obtains
for instance from Corollary~\ref{co:adjoint} by embedding $\mathbf
K^b(\C)/{\langle\Cone\Phi\rangle}$ into the derived category $\mathbf
D(R)$.

\begin{lem}\label{le:universal}
The composite $Q\colon\mathbf K^b(\C)\to\mathbf
K^b(\C)/{\langle\Cone\Phi\rangle}\to\T$ has the following properties.
\begin{enumerate}
\item $\T(\Si^n (QX),QY)=0$ for all
$X,Y$ in $\C$ and $n>0$.
\item The functor $Q|_\C\colon\C\to\T$ factors through the
localization $\C\to\C[\Phi^{-1}]$.
\item The functor $\C[\Phi^{-1}]\to\T$ extends to an exact functor $\mathbf
K^b(\C[\Phi^{-1}])\to\T$.
\item The functor $-\otimes_RS\colon\mathbf K^b(\C)\to\mathbf
K^b(\C[\Phi^{-1}])$ factors through $Q\colon\mathbf K^b(\C)\to\T$
\end{enumerate}
\end{lem}
\begin{proof}
(1) The functor $Q\colon\mathbf K^b(\C)\to\T$ is a cohomological
quotient functor.  Thus we can apply Corollary~\ref{co:adjoint} and
obtain a fully faithful and exact `right adjoint' $Q'\colon
\T\to\mathbf D(R)$ such that
$$\T(QA,B)\cong\mathbf D(R)(A,Q'B)\quad\textrm{for all $A\in\mathbf K^b(R)$ and $B\in\T$}.$$
To compute $\T(\Si^n (QX),QY)$ for $X,Y$ in $\C$, it is sufficient to consider the case $X=R=Y$.
We have
$$\T(\Si^n (QR),QR)\cong\mathbf D(R)(\Si^n R,(Q'\comp Q)R)\cong
H^{-n}((Q'\comp Q)R).$$ Now we apply Corollary~3.31 from \cite{NR}
which says that $\Tor_n^R(S,S)=0$ for all $n>0$ implies
$H^{-n}((Q'\comp Q)R)=0$ for all $n>0$.

(2) The functor $\mathbf K^b(\C)\to\mathbf
K^b(\C)/{\langle\Cone\Phi\rangle}$ makes the maps in $\Phi$ invertible
by sending the objects in $\Cone\Phi$ to zero. Therefore $\C\to\T$
factors through the localization $\C\to\C[\Phi^{-1}]$.

(3) This follows from the `universal property' of the homotopy
category $\mathbf K^b(\C[\Phi^{-1}])$ which is the main result in
\cite{KV}. More precisely, any additive functor $F\colon\D\to\bar \A$
from an additive category $\D$ to the stable category of a Frobenius
category $\A$ extends to an exact functor $\mathbf K^b(\D)\to\bar\A$
provided that $\bar \A(\Si^n(FX),FY)=0$ for all $X,Y$ in $\D$ and
$n>0$. Note that we are using (1) and the fact that $\T$ can be
embedded into the stable category of a Frobenius category.

(4) We have $\mathbf K^b(S)=\mathbf K^b(\C[\Phi^{-1}])$ and
$X\otimes_RS=0$ for all $X\in\Cone\Phi$ since $R\to S$ is the universal
localization with respect to $\Phi$. Thus $-\otimes_RS$ factors
through the quotient functor $\mathbf K^b(\C)\to\mathbf
K^b(\C)/{\langle\Cone\Phi\rangle}$. Moreover, $-\otimes_RS$ factors
through $\T$ since idempotents in $\mathbf K^b(\C[\Phi^{-1}])$ split.
\end{proof}

The following commutative diagram summarizes our findings from
Lemma~\ref{le:universal}.

\begin{equation}\label{eq:FG}
\xymatrix{\C\ar[d]\ar[r]&\C[\Phi^{-1}]\ar[r]\ar[d]&\mathbf
K^b(\C[\Phi^{-1}])\ar[ld]^F\\ \mathbf
K^b(\C)\ar[d]^T\ar[r]^Q&\T\ar[ld]^G\\ \mathbf K^b(\C[\Phi^{-1}])}
\end{equation}

\begin{proof}[Proof of Theorem~\ref{th:universal}] 
We want to apply Theorem~\ref{th:maplifting} and use the diagram
(\ref{eq:FG}). More precisely, we need to show that $\Ann T=\Ann Q$,
because this implies condition (2) in Theorem~\ref{th:maplifting}
since $T=-\otimes_RS$ and $\Ann Q$ is generated by identity maps.  We
claim that $Q=F\comp T$.  This follows from the `universal property'
of the homotopy category $\mathbf K^b(\C)$ since $Q|_\C=F\comp T|_\C$;
see \cite{Ke2}. We obtain that $Q= F\comp G\comp Q$, and this implies
$F\comp G\cong\Id_\T$ since both functors agree on $\mathbf
K^b(\C)/{\langle\Cone\Phi\rangle}$.  Thus $G$ is faithful and we
conclude that $\Ann T=\Ann Q$. This completes the proof.
\end{proof}

\begin{rem}
I conjecture that Theorem~\ref{th:universal} remains true if one
replaces `universal localization' by `homological localization'.
\end{rem}

\section{Almost derived categories}

Almost rings and modules have been introduced by Gabber and Ramero
\cite{GR}.  Here, we analyze their formal properties and introduce
their analogue for derived categories. Let us start with a piece of notation.
Given a class $\Phi$ of maps in some additive category $\C$, we denote by
$$\Phi^\perp=\{X\in\C\mid \mbox{$\C(\p,X)=0$ for all $\p\in\Phi$}\}$$
the full subcategory of objects which are {\em annihilated} by $\Phi$.

Throughout this section we fix an associative ring $R$. We view
elements of $R$ as maps $R\to R$. Thus $\mathfrak a^\perp$ for any
ideal $\mathfrak a$ of $R$ denotes the category of $R$ modules
which are annihilated by $R$.

The formal essence of an almost module category can be formulated as
follows.

\begin{prop}\label{pr:almost}
Let $\A$ be a full subcategory of a module category $\Mod R$.  Then
the following are equivalent.
\begin{enumerate}
\item $\A$ is a Serre subcategory, and the inclusion has a left and a right adjoint.
\item There exists an idempotent ideal $\mathfrak a$ of $R$ such that
$\A=\mathfrak a^\perp$.
\end{enumerate}
In this case, the quotient category $\Mod R/\A$ is the category of
almost modules with respect to $\mathfrak a$, which is denoted by
$\Mod (R,\mathfrak a)$.
\end{prop}
\begin{proof} 
The proof of the first part is straightforward; see for instance
\cite[Proposition~7.1]{A}.  The second part is just the definition of
an almost module category from \cite{GR}.
\end{proof}

The following result is the analogue of Proposition~\ref{pr:almost}
for triangulated categories.

\begin{thm}\label{th:almost}
Let $\R$ be a full subcategory of a compactly generated triangulated
category $\S$. Then the following are equivalent.
\begin{enumerate}
\item $\R$ is a triangulated subcategory, and the inclusion has a left
and a right adjoint.
\item There exists an idempotent ideal $\mathfrak I$ of $\S_c$
satisfying $\Si\mathfrak I=\mathfrak I$, such that $\R=\mathfrak I^\perp$.
\end{enumerate}
In this case, the left adjoint of the quotient functor
$\S\to\S/\R$ identifies $\S/\R$ with a smashing subcategory of $\S$.
Moreover, every smashing subcategory of $\S$ arises in this way.
\end{thm}
\begin{proof}
Let us denote by $F\colon \R\to\S$ the inclusion functor.

(1) $\Rightarrow$ (2): The left adjoint $E\colon \S\to \R$ of the
inclusion $\R\to\S$ is a smashing localization functor since
$E\comp F\cong\Id_\R$. It follows from Theorem~\ref{th:smash2} and its
Corollary~\ref{co:exact} that $\mathfrak I=\S_c\cap\Ann E$ is an
idempotent ideal satisfying $\Si\mathfrak I=\mathfrak I$ and
$\mathfrak I^\perp=\R$.

(2) $\Rightarrow$ (1): Lemma~\ref{le:perp} implies that $\mathfrak
I^\perp=\R$ is a triangulated subcategory. Let us replace $\mathfrak
I$ by the ideal $\mathfrak J$ of all maps in $\S_c$ annihilating
$\R$. Thus $\mathfrak J$ is a cohomological ideal satisfying
$\Si\mathfrak J=\mathfrak J$ and $\mathfrak J^\perp=\R$. The proof of
Theorem~\ref{th:smash2} shows that $\R$ is perfectly cogenerated. Thus
the inclusion $F\colon\R\to\S$ has a left adjoint by
Corollary~\ref{co:brown}, since $F$ preserves all
products. Theorem~\ref{th:smash1} implies that $\R$ is compactly
generated. Thus $F$ has a right adjoint by the dual of
Corollary~\ref{co:brown}, since $F$ preserves all coproducts.

Now let us prove the second part. Suppose first that (1) -- (2) hold.
Then the left adjoint $E\colon \S\to\R$ of the inclusion $\R\to\S$ is
a smashing localization functor. It follows that the left adjoint of
the quotient functor $\S\to \S/\R$ identifies $\S/\R$ with $\Ker E$,
which is by definition a smashing subcategory.

Finally, suppose that $\T$ is a smashing subcategory of $\S$.  Let
$\mathfrak I$ be the idempotent ideal of all maps in $\S_c$ which
factor through some object in $\T$. Then $\R=\mathfrak I^\perp$ is a
triangulated subcategory of $\S$, and the inclusion $\R\to\S$ has a
left and a right adjoint. It follows that the left adjoint of the
quotient functor $\S\to \S/\R$ identifies $\S/\R$ with $\T$.  This
finishes the proof.
\end{proof}

Let us complete the parallel between module categories and derived
categories.  Thus we consider the unbounded derived category $\mathbf
D(R)$ of the module category $\Mod R$.  Comparing the statements of
Proposition~\ref{pr:almost} and Theorem~\ref{th:almost}, we see that
the formal analogue of an almost module category is a triangulated
category of the form
$$\mathbf D(R,\mathfrak I)=\mathbf D(R)/(\mathfrak I^\perp)$$ for some
idempotent ideal $\mathfrak I$ of $\mathbf K^b(R)$ satisfying
$\Si\mathfrak I=\mathfrak I$.  We call such a category an {\em almost
derived category}.

Next we show that the derived category of an almost module
category is an almost derived category.

\begin{cor}\label{co:almost}
Let $R$ be a ring and $\mathfrak a$ be an idempotent ideal such that
$\mathfrak a\otimes_R\mathfrak a$ is flat as left $R$-module. Denote
by $\mathfrak A$ the maps in $\mathbf K^b(R)$ which annihilate all
suspensions of the mapping cone of the natural map $\mathfrak
a\otimes_R\mathfrak a\to R$.  Then we have
$$\mathfrak A^2=\mathfrak A\quad\textrm{and}\quad\mathbf
D(\Mod(R,\mathfrak a))=\mathbf D(R,\mathfrak A).$$
\end{cor}
\begin{proof}
The quotient functor $F\colon\Mod R\to\Mod (R,\mathfrak
a)$ has a left adjoint $E$, and we have
$$(E\comp F)M=M\otimes_R(\mathfrak a\otimes_R\mathfrak a);$$ see for
instance \cite[p.\ 200]{St}. The extra assumption on $\mathfrak a$
implies that $E$ is exact.  Taking derived functors, we obtain an
adjoint pair of exact functors
$$\mathbf RF\colon \mathbf D(\Mod R)\lto\mathbf D(\Mod (R,\mathfrak
a))\quad\textrm{and}\quad \mathbf LE\colon \mathbf D(\Mod
(R,\mathfrak a))\lto\mathbf D(\Mod R)$$ such that $\mathbf RF\comp
\mathbf LE\cong\Id_{\mathbf D(\Mod (R,\mathfrak a))}$.  It
follows that $\mathbf LE$ identifies $\mathbf D(\Mod
(R,\mathfrak a))$ with a smashing subcategory $\R$ of $\mathbf D(R)$.
Now observe that the mapping cone $\Cone (\mathfrak
a\otimes_R\mathfrak a\to R)$ generates $\mathbf D(R)/\R$. In fact, the
canonical map $$R\lto\Cone (\mathfrak a\otimes_R\mathfrak a\to R)$$ is
an isomorphism in $\mathbf D(R)/\R$ since $(\mathbf LE\comp
\mathbf RF)R=\mathfrak a\otimes_R\mathfrak a$. Therefore $\mathfrak A$
is the exact ideal corresponding to $\R$ which is idempotent by
Corollary~\ref{co:exact}.  Moreover, the localization functor $\mathbf
RF$ identifies $\mathbf D(R)/(\mathfrak A^\perp)$ with $\mathbf
D(\Mod(R,\mathfrak a))$.
\end{proof}

We know that the derived category of a ring is compactly generated.
This is no longer true for almost derived categories \cite{Ke}.  In
deed, the telescope conjecture expresses the fact that all almost
derived categories are compactly generated.

\begin{cor}
The telescope conjecture holds for the derived category $\mathbf D(R)$
of a ring $R$ if and only if every almost derived category $\mathbf
D(R,\mathfrak I)$ is a compactly generated triangulated category.
\end{cor}

\begin{appendix}
\section{Epimorphisms of additive categories}\label{ap:epi}

An additive functor $F\colon\C\to\D$ between additive categories is
called an {\em epimorphism of additive categories}, or simply an {\em
epimorphism}, if $G\comp F=G'\comp F$ implies $G=G'$ for any pair
$G,G'\colon\D\to \E$ of additive functors. In this section we
characterize epimorphisms of additive categories in terms of functors
between their module categories.  This material is classical \cite{M},
but we need it in a form which slightly generalizes the usual
approach.

\begin{lem}\label{le:epi1a}
  Let $F\colon\C\to\D$ be an additive functor between additive
  categories.  Suppose the restriction $F_*\colon\Mod\D\to\Mod\C$ is
  full and $F$ is surjective on objects. Then $F$ is an epimorphism.
\end{lem}
\begin{proof} 
Let $G,G'\colon\D\to \E$ be a pair of additive
functors satisfying $G\comp F=G'\comp F$. Clearly, $G$ and $G'$
coincide on objects since $F$ is surjective on objects. Now choose a
map $\a\colon X\to Y$ in $\D$. We need to show that $G\a=G'\a$. The
functor $G'$ induces a $\C$-linear map
$$\g\colon F_*\D(-,Y)\lto (F_*\comp G_*)\E(-,GY),$$
which is defined
by $$\g_C\colon\D(FC,Y)\lto\E(G(FC),GY),\quad \p\mapsto G'\p$$
for
each $C$ in $\C$.  The fact that $F_*$ is full implies that $\g=F_*\d$
for some $\D$-linear map $\d\colon \D(-,Y)\to G_*\E(-,GY)$. In
particular, $\d_X=\g_C$ for some $C$ in $\C$ satisfying $FC=X$. Thus
we obtain the following commutative diagram
$$\xymatrix{
\D(Y,Y)\ar[rr]^{\d_{Y}}\ar[d]^{\D(\a,Y)}&&\E(GY,GY)\ar[d]^{\E(G\a,GY)}\\
\D(X,Y)\ar[rr]^{\d_{X}}&&\E(GX,GY)\\ }$$ which shows
$G\a=G'\a$ if we apply it to $\id_Y$. We conclude that $G=G'$.
\end{proof}

\begin{lem}\label{le:epi1b}
  Let $F\colon\C\to\D$ be an additive functor between additive
  categories. Suppose $F$ is an epimorphism and bijective on objects.  Then the
 restriction $F_*\colon\Mod\D\to\Mod\C$ is fully faithful.
\end{lem}
\begin{proof} 
Let $M,N$ be a pair of $\D$-modules. We need to
show that the canonical map
$$(F_*)_{M,N}\colon\Hom_\D(M,N)\lto\Hom_\C(F_*M,F_*N)$$ is
bijective. Given a family $\p=(\p_X)_{X\in\D}$ of maps $\p_X\colon
MX\to NX$, we define a $\D$-module $H_\p$ by $$H_\p X=MX\amalg NX
\quad\textrm{and}\quad H_\p\a=\left[ \begin{matrix}M\a&0\\
N\a\comp\p_Y-\p_X\comp M\a&N\a\end{matrix}\right]$$ for each object
$X$ and each map $\a\colon X\to Y$ in $\D$.  Note that
$$(\p_X)_{X\in\D}\colon M\lto N$$ is $\D$-linear if and only if
$H_\p=M\amalg N$.

To prove that $(F_*)_{M,N}$ is surjective, fix a $\C$-linear map
$$\psi=(\psi_{X})_{X\in\C}\colon F_*M\lto F_*N.$$
For each $X$ in $\D$
put $\p_X=\psi_{F^{-1}X}$.  We have $H_\p\comp F=(M\amalg N)\comp F$ since
$\psi$ is $\C$-linear.  Thus $\p$ is $\D$-linear because $H_\p\comp
F=(M\amalg N)\comp F$ implies $H_\p=M\amalg N$. We have $F_*\p=\psi$
and conclude that the map $(F_*)_{M,N}$ is surjective.

To prove that $(F_*)_{M,N}$ is injective, choose a non-zero map
$\p\colon M\to N$.  Thus $\Im\p\neq 0$. We have
$\Im(F_*\p)=F_*(\Im\p)\neq 0$ and therefore $F_*\p\neq 0$.  It follows
that $(F_*)_{M,N}$ is injective.
\end{proof}

\begin{lem}\label{le:epi2}
  Let $F\colon\C\to\D$ be an additive functor between additive
  categories.  If $F$ is an epimorphism, then $F$ is surjective on
  objects.
\end{lem}
\begin{proof} 
Suppose there is an object $D$ in $\D$ which does not belong to the
image of $F$. We construct a new additive category $\E$ which contains
$\D$ as a full subcategory and has one additional object, denoted by
$D'$.  Let $\E(X,D')=\D(X,D)$ and $\E(D',X)=\D(D,X)$ for all $X$ in
$\D$, and let $\E(D',D')=\D(D,D)$. Now define $G\colon\D\to \E$ to be
the inclusion, and define $G'\colon\D\to \E$ by $G'X=GX$ for all $X$
in $\D$, except for $X=D$, where we put $G'D=D'$. Clearly, $G\comp
F=G'\comp F$ but $G\neq G'$. Thus an epimorphism is surjective on
objects.
\end{proof}

\begin{lem}\label{le:epi3}
  Let $F\colon\C\to\D$ be an additive functor between small additive
  categories. If the restriction $F_*\colon\Mod\D\to\Mod\C$ is
  faithful, then every object in $\D$ is a direct factor of some
  object in the image of $F$.
\end{lem}
\begin{proof}
  The restriction $F_*$ has a left adjoint
  $F^*\colon\Mod\C\to\Mod\D$. The assumption on $F_*$ implies that for
  each $\D$-module $M$ the natural map $(F^*\comp F_*)M\to M$ is an
  epimorphism. Now fix an object $Y$ in $\D$.  Every module is a
  quotient of a coproduct of representable functors.  Thus we have an
  epimorphism
$$\coprod_{i\in\La}\C(-,X_i)\lto F_*\D(-,Y),$$
and applying $F^*$ induces an epimorphism
$$\coprod_{i\in\La}\D(-,FX_i)\lto(F^*\comp F_*)\D(-,Y)\lto\D(-,Y).$$ 
Using Yoneda's lemma, we
see that $Y$ is a direct factor of $F(\coprod_{i\in\Ga}X_i)$ for some
finite subset $\Ga\subseteq\La$.
\end{proof}

\begin{prop}\label{pr:epi4}
  Let $F\colon\C\to\D$ be an additive functor between small additive
  categories.  Then $F_*\colon\Mod\D\to\Mod\C$ is fully faithful if and
  only if there is a factorization $F=F_2\comp F_1$ such that 
\begin{enumerate}
\item $F_1$ is an epimorphism and bijective on objects, and
\item $F_2$ is fully faithful and every object in $\D$ is a direct
factor of some object in the image of $F_2$.
\end{enumerate}
\end{prop}
\begin{proof} 
Suppose first that $F_*$ is fully faithful.
We define a factorization
$$\C\stackrel{F_1}\lto\D'\stackrel{F_2}\lto\D$$ as follows.  The
objects of $\D'$ are those of $\C$ and $F_1$ is the identity on
objects. Let $$\D'(X,Y)=\D(FX,FY)$$ for all $X,Y$ in $\C$, and let
$F_1\a=F\a$ for each map $\a$ in $\C$. The functor $F_2$ equals $F$ on
objects and is the identity on maps. It follows that $F_2$ is fully
faithful and surjective up to direct factors on objects, by
Lemma~\ref{le:epi3}. Thus $(F_2)_*$ is fully faithful, and this
implies that $(F_1)_*$ is fully faithful, since
$F_*=(F_2)_*\comp(F_1)_*$. We conclude from Lemma~\ref{le:epi1a} that
$F_1$ is an epimorphism.

Now suppose $F$ admits a factorization $F=F_2\comp F_1$ satisfying (1)
and (2). Then $(F_1)_*$ is fully faithful by Lemma~\ref{le:epi1b}, and
$(F_2)_*$ is automatically fully faithful. Thus $F_*$ is fully faithful.
\end{proof}

The property of being an epimorphism is invariant under enlarging the
universe.

\begin{lem}\label{le:universe}
Let $\mathfrak U$ and $\mathfrak V$ be universes in the sense of
Grothendieck \cite[I.1]{G}, and suppose $\mathfrak U\subseteq\mathfrak
V$. If $F\colon\C\to\D$ is an epimorphism of additive $\mathfrak
U$-categories, then $F$ is an epimorphism of additive $\mathfrak
V$-categories.
\end{lem}
\begin{proof}
Let $G,G'\colon\D\to \E$ be a pair of additive functors into a
$\mathfrak V$-category $\E$ satisfying $G\comp F=G'\comp F$. We denote
by $\F$ the smallest additive subcategory of $\E$ containing the image
of $G$ and $G'$. Observe that $\F$ is a $\mathfrak U$-category since
$\D$ is a $\mathfrak U$-category. Thus the restrictions $\D\to \F$ of
$G$ and $G'$ agree by our assumption on $F$. It follows that $G=G'$.
\end{proof}

\section{The abelianization of a triangulated category}\label{ap:abel}

Let $\C$ be a triangulated category. In this section we discuss some
properties of the {\em abelianization} $\mod\C$ of $\C$. Most of this
material can be found in work of Freyd \cite{F} about the
formal properties of the stable homotopy category.

\begin{lem}
  Let $\C$ be a triangulated category. Then the category $\mod\C$ is
  an abelian {\em Frobenius category}, that is, there are enough
  projectives and enough injectives, and both coincide.
\end{lem}
\begin{proof}
The representable functors are projective objects in $\mod\C$ by
Yoneda's lemma. Thus $\mod\C$ has enough projectives.  Using the fact
that the Yoneda functors $\C\to\mod\C$ and $\C^\op\to\mod(\C^\op)$ are
universal cohomological functors, we obtain an equivalence
$(\mod\C)^\op\to\mod(\C^\op)$ which sends $\C(-,X)$ to $\C(X,-)$ for
all $X$ in $\C$. Thus the representable functors are injective
objects, and $\mod\C$ has enough injectives.
\end{proof}

The triangulated structure of $\C$ induces some additional structure
on $\mod\C$. This involves the equivalence
$\Si^\star\colon\mod\C\to\mod\C$ which extends $\Si\colon\C\to\C$.
By abuse of notation, we identify $\Si^\star=\Si$.  Using this
internal grading, the category $\mod\C$ is $(3,-1)$-periodic \cite{F}.
Thus we obtain a canonical extension $\k_M$ in $\Ext_\C^3(\Si M,M)$
for every module $M$ in $\mod\C$. Under some additional assumptions,
this extension is induced by a Hochschild cocycle of degree $(3,-1)$;
it plays a crucial role in \cite{BKS}.

\begin{prop}\label{pr:kappa}
Let $\C$ be a triangulated category.
\begin{enumerate}
\item Given a pair $M,N$ of objects in $\mod\C$, there is a natural
map
$$\k_{M,N}\colon\Hom_\C(M,N)\lto\Ext^3_\C(\Si M,N)$$ and we write
$\k_N=\k_{N,N}(\id_N)$. 
\item Let
  $$\Delta\colon X\stackrel{\a}\lto Y\stackrel{\b}\lto
  Z\stackrel{\g}\lto \Si X$$
  be a sequence of maps in $\C$ and let
  $N=\Ker\C(-,\a)$.  Then $\Delta$ is an exact triangle if and only if
  the map $\g$ induces a map $\C(-,Z)\lto \Si N$ such that the
  sequence
$$0\lto N\lto \C(-,X)\xrightarrow{\C(-,\a)}
\C(-,Y)\xrightarrow{\C(-,\b)} \C(-,Z)\lto \Si N\lto 0$$
is exact in $\mod\C$ and represents $\k_N$.
\end{enumerate}
\end{prop}
\begin{proof}
(1) Let $M=\Coker\C(-,\b)$ be a an object in $\mod\C$ and complete
$\b\colon Y\to Z$ to an exact triangle $X\to Y\to Z\to\Si X$ to obtain
a projective resolution
$$\cdots\to\C(-,Y)\to\C(-,Z)\to \C(-,\Si X)\to \C(-,\Si Y)\to\C(-,\Si
Z)\to \Si M\to 0$$ of $\Si M$.  The map $\k_{M,N}$ takes by definition
a map $\p\colon M\to N$ to the element in $\Ext^3_\C(\Si M,N)$ which
is represented the composition of $\p$ with the projection $\C(-,Z)\to
M$.

(2) Fix a sequence
$$\Delta\colon X\stackrel{\a}\lto Y\stackrel{\b}\lto
Z\stackrel{\g}\lto \Si X$$
in $\C$ and let $N=\Ker\C(-,\a)$. Suppose
first that $\Delta$ is an exact triangle. The definition of $\k_N$
implies that the induced sequence
$$\e_\Delta\colon 0\lto N\lto \C(-,X)\xrightarrow{\C(-,\a)}
\C(-,Y)\xrightarrow{\C(-,\b)} \C(-,Z)\lto \Si N\lto 0$$
is exact in
$\mod\C$ and represents $\k_N$. Conversely, suppose that $\e_\Delta$
is exact and represents $\k_N$. Complete $\a$ to an exact triangle
$$\Delta'\colon X\stackrel{\a}\lto Y\stackrel{\b'}\lto
Z'\stackrel{\g'}\lto \Si X$$ in $\C$. We use dimension shift and replace
both sequences $\e_\Delta$ and  $\e_{\Delta'}$ by  short exact sequences
$$0\to\Omega^{-2}N \to \C(-,Z)\to \Si N\to 0\quad\textrm{and}\quad
0\to\Omega^{-2}N \to \C(-,Z')\to \Si N\to 0.$$ These represent the same
element in $\Ext_\C^1(\Si N,\Omega^{-2}N)$, and we obtain therefore an
isomorphism $\p\colon Z\to Z'$ which induces an isomorphism of
triangles $\Delta\to\Delta'$. Thus $\Delta$ is an exact triangle.
\end{proof}

Let us explain a more conceptual way to understand the natural map
$\k_{M,N}$. To this end denote by $\umod\C$ the {\em stable category} of
$\mod\C$, that is, the objects are those of $\mod\C$ and 
$$\uHom_\C(M,N)=\Hom_\C(M,N)/{\mathfrak P(M,N)}$$ where $\mathfrak P$
denotes the ideal of all maps in $\mod\C$ which factor through some
projective object. Taking syzygies in $\mod\C$ induces an equivalence
$$\Omega\colon\umod\C\lto\umod\C$$ since $\mod\C$ is a Frobenius
category.  Moreover, $$\uHom_\C(\Omega^n
M,N)\cong\Ext_\C^n(M,N)\cong\uHom_\C(M,\Omega^{-n}N) \quad\textrm{for
all $M,N$ and $n>0$}.$$ The map $\k_{M,N}$ induces a natural isomorphism
$$\uHom_\C(M,N)\lto\Ext^3_\C(\Si M,N),$$ and composing this with the
natural isomorphism
$$\Ext_\C^3(\Si M,N)\lto\uHom_\C(\Si M,\Omega^{-3}N)$$
induces a natural isomorphism between
$$\bar\Si\colon\umod\C\lto\umod\C
\quad\textrm{and}\quad\Omega^{-3}\colon\umod\C\lto\umod\C.$$ Note that
the natural map $\k_{M,N}$ can be reconstructed from the natural
isomorphism $\bar\Si\cong\Omega^{-3}$.

\end{appendix}

\end{document}